\documentclass[twoside,12pt]{article}
\usepackage{amsmath}
\usepackage{amssymb}
\usepackage{mathrsfs}
\usepackage{amsthm}
\usepackage{graphicx}

\def\fukao{\color{red}}
\let\fukao\relax

\def\risei{\color{red}}
\let\risei\relax

\title{Time-dependence of the threshold function in the perfect plasticity model}

\author{
Takeshi Fukao\\
Department of Mathematics, Faculty of Education\\
Kyoto University of Education\\
1~Fujinomori, Fukakusa, Fushimi-ku, Kyoto~612-8522 Japan\\
E-mail: \texttt{fukao@kyokyo-u.ac.jp}
\and \\ 
Risei Kano\\
Research and Education Faculty, Humanities and Social Science Cluster\\
Kochi University\\
2-5-1~Akebono-cho, Kochi~780-8520 Japan\\
E-mail: \texttt{kano@kochi-u.ac.jp}}
 
\date{}

\newcommand\testopari{\sc Takeshi Fukao and Risei Kano}
\newcommand\testodispari{\sc Perfect plasticity model with time dependent constraints}
\begin{document}

\maketitle

\begin{abstract}

This paper discusses the time-dependence of the threshold function in 
the perfect plasticity model. 
In physical terms, it is natural that 
the threshold function depends on some unknown variable. 
Therefore, it is meaningful to discuss 
the well-posedness of this function under the weaker assumption of 
time-dependence. 
Time-dependence is also interesting from the viewpoint of the abstract evolution equation. 
To prove the existence of a solution to the perfect plasticity model, 
the recent abstract theory under the continuous class with respect to time 
is used. \\

\vspace{2mm}

\vspace{2mm}
\noindent \textbf{AMS (MOS) subject clas\-si\-fi\-ca\-tion:} 35K61, 35K25, 35D30, 80A22.\\

\end{abstract}

\newpage

\section{Introduction}
\setcounter{equation}{0}

{\fukao In 2004, U.\ {S}tefanelli advanced an interesting subject \cite{Ste04} with respect to 
the abstract theory of evolution inclusion, 
more precisely the relationship between \emph{nonlocal quasivariational problem} and 
\emph{hardening problem} on the elastoplastic materials. 
It is well known that the essence of this physical problem is characterized 
by the convex constraint for unknown stress $\sigma$, mathematically. 
His idea was the characterization of the constraint set depends on the unknown stress itself, nonlocally. 
In this paper, in order to progress this direction from the view point of abstract theory of evolution equation, 
we focus on the time-dependent constraint set, more precisely 
we consider the prototype model of the plasticity  
namely the \emph{perfect plasticity model} 
of G.\ {D}uvaut and J.\ L.\ {L}ions \cite{DL76}, and focus on the \emph{time-dependence} 
of the constraint set of the 
following form: 
\begin{equation*}
	\dot{\varepsilon}= \dot{\sigma} + \lambda, \quad  \lambda \in \partial I_{Z}(\sigma),
\end{equation*}
where $\varepsilon :=\varepsilon(u)$ stands for the strain tensor and the symbol 
$\dot{\varepsilon }$ and $\dot{\sigma }$ mean the time derivatives of 
$\varepsilon $ and $\sigma $, respectively. 
We chose the suitable set $Z$ with related to {V}on {M}ises model or 
{T}resca model and so on. 
The prototype model in \cite{DL76}, the set $Z$ is defined by 
\begin{equation*} 
	Z:=
	\left\{ 
	\tau \in \mathbb{R}^{3\times3}_{\textrm{sym}} : 
	\frac{1}{2} \bigl| \tau^D\bigr|^2 \le g  \right\},
\end{equation*}
where $\mathbb{R}^{3\times3}_{\textrm{sym}}$ stands for the $3 \times 3$ symmetric matrix, 
$\tau ^D$ stands for the deviator of $\sigma $ defined by
$\tau^D_{ij}:=
\tau_{ij}-(1/3)\sum_{k=1}^3\tau_{kk}
\delta_{ij}$ for $i,j=1,2,3$
and 
$g$ is a positive constant that represents the threshold value of the 
elasto-plasticity. 
In this paper we focus on the the case that $g$ is a function and find the time-space regularity assumption 
for solving the perfect plastic model, based on the abstract theory of time-dependent evolution inclusion. 
Indeed, we can find very interesting model \cite{BBW15}, 
where $g$ depends on some another unknown temperature. In this sense, the time-dependent model 
is a first step for this direction. } 

The time-dependent evolution inclusion is an interesting subject in the 
theory of the abstract evolution equations. Here, 
\emph{time-dependence} includes the time-dependent effective domain, i.e.\ we consider the 
following evolution inclusion:
\begin{equation}  
	u'(t) + \partial \varphi ^t \bigl( u(t) \bigr) \ni f(t) 
	\quad {\rm in}~H, \ {\rm for~a.a.~} t \in (0,T), 
	\quad u(0)=u_0 \quad {\rm in}~H,
	\label{ee1}
\end{equation} 
 where $H$ is a real {H}ilbert space and $T>0$ is a finite time; 
$\partial \varphi^t$ is a subdifferential of some proper, lower semi-continuous, convex 
functional $\varphi ^t$ on $H$, which depends on time $t \in [0,T]$. 
More precisely, the effective domain 
$D(\varphi ^t)\subset H$ moves depending on time. 
The dynamics have some constraints that depend on time. 
The well-posedness of the above evolution inclusion has been 
studied in many papers \cite{Bir74, Bir74b, Bre72, FK13, Ken75, Ken81, Yam76} under 
various settings of the time-dependence. In this paper, 
we recall one result given by {F}ukao--{K}enmochi \cite{FK13} 
in terms of a weakly time-dependent constraint. 
Moreover, we will apply this result to a version of 
the classical perfect plasticity problem \cite{DL76}. 
The {\fukao motivation} of this study is to apply the theory of 
quasivariational inequality (see, e.g.\ {\risei{\cite{KMK09, KL02, KM98, Ste04}}}) to 
the hardening problem, {\fukao we will treat this subject in the forthcoming paper.}
In this point, the special case of \eqref{ee1} is very important, 
namely $\varphi ^t:=I_{K(t)}$, where 
$I_{K(t)}$ is the indicator function of some closed convex set $K(t)$, 
which is called the \emph{{M}oreau sweeping process}. 
In the 1970s, J.\ J.\ {M}oreau studied a class of evolution inclusion of the form
\begin{equation*}  
	u'(t) + \partial I_{K(t)} \bigl( u(t) \bigr) \ni 0 
	\quad {\rm in}~H, \ {\rm for~a.a.~} t \in (0,T), 
	\quad u(0)=u_0 \quad {\rm in}~H,
\end{equation*} 
(see, e.g.\ \cite{Mor71, Mor77}). 
This problem has been studied in various scenarios (see, e.g.\ 
\cite{BKS04, KL02, KL09, KR11, KM98, Mon93, Rec11, Rec15, Vla13}).

In Section~2, we introduce the original problem of 
perfect plasticity (see, e.g.\ \cite{DL76, Kre96, Sho97}). 
After setting function spaces with some useful properties, we define our solution for 
three cases. Moreover, we give three theorems as examples of the abstract theory of 
evolution inclusion.

In Section~3, to prove the main theorem, we recall the abstract theory of 
evolution equations with time-dependence. 
The first theorem (Theorem 2.1) is a consequence of the well-known theory. 
Moreover, using this result, we can discuss the well-posedness of
the weak variational inequality (Theorem 2.2). In both cases, 
we need to consider approximations under two parameters, $\kappa, \nu >0$.

In Section~4, we consider the problem without the approximate parameter $\kappa >0$. 
In this case, the problem is the same as the well-known {M}oreau sweeping process. 
{\fukao Therefore, we do not consider the limiting procedure $\kappa \to 0$ from previous section, 
we directly treat the problem as the {M}oreau sweeping process. }

A detailed index of sections and subsections is as follows:

\begin{itemize}
 \item[1.] Introduction
 \item[2.] Main theorems
	\begin{itemize}
	 \item[2.1.] Original problem
	 \item[2.2.] Notation
	 \item[2.3.] Definition of the solution and theorem
	\end{itemize}
 \item[3.] Proof of the main theorems
	\begin{itemize}
	 \item[3.1] Evolution inclusion with time-dependent domain
	 \item[3.2] Proof of Theorem 2.1
	 \item[3.3] Proof of Theorem 2.2
	\end{itemize}
 \item[4.] {\fukao Viscous} perfect plasticity model
	\begin{itemize}
	 \item[4.1] Auxiliary problems 
	 \item[4.2] Proof of Theorem 2.3
	\end{itemize}
\end{itemize}

\section{Main theorems}
\setcounter{equation}{0}

In this paper, we consider the time-dependence of the 
constraint related to the perfect plasticity model. 
In this section, we state the main results. 
{\fukao Throughtout of this paper we use the notation 
$\boldsymbol{u}:=(u_1,u_2,u_3)$ for the vector and  
$\tau :=\{ \tau _{ij} \}$ for the tensor with the following matrix representation
\begin{equation*}
	\{ \tau _{ij} \}=
	\begin{pmatrix}
	\tau _{11}	& \tau _{12}	& \tau _{13}	\\
	\tau _{21}	& \tau _{22}	& \tau _{23}	\\
	\tau _{31}	& \tau _{32}	& \tau _{33}
	\end{pmatrix}
\end{equation*}
with $\tau _{i\cdot }:=(\tau _{i1}, \tau_{i2}, \tau _{i3})$ for $i=1,2,3$. 
}

\subsection{Original problem}

In this subsection, we recall the well-known classical problem of 
perfect plasticity \cite{DL76}. 
The unknown functions 
$\boldsymbol{u}=\boldsymbol{u}(t,x)$ 
and $\sigma=\sigma(t,x)$ describe the displacement and stress, respectively, 
in the interior of a medium that occupies a smooth region $\Omega \subset \mathbb{R}^3$. 
The boundary $\Gamma :=\partial \Omega$ consists of $\Gamma =\Gamma _D \cup\Gamma _N$, 
where $\Gamma _D \cap \Gamma _N = \emptyset$ with 
{\fukao measures} 
$|\Gamma _D|>0$ and $|\Gamma _N|>0$. Moreover, 
$\boldsymbol{\nu }$ denotes the unit normal vector outward from $\Gamma $. 
The function $\varepsilon(\boldsymbol{u}):=\{ \varepsilon _{ij}\}$ represents
the strain with respect to deformation, defined by
$\varepsilon _{ij}:=(1/2)((\partial u_i/\partial x_j)
+(\partial u_j/\partial x_i))$ for $i,j=1,2,3$. 
We wish to find $\boldsymbol{v}:=\partial \boldsymbol{u}/\partial t$ 
and $\sigma $ satisfying 
\begin{align}
	&\frac{\partial \boldsymbol{v} }{\partial t} 
	= \mathbf{div}\sigma + \boldsymbol{f} \quad \mbox{in}~Q:=(0,T) \times \Omega \label{dleqv},\\
	&\frac{\partial \sigma}{\partial t}+\partial I_{ Z}(\sigma+\sigma_{*})\ni \varepsilon 
	 (\boldsymbol{v})+h \quad \mbox{in}~Q,\label{dleqs}
\end{align}
under some initial and boundary conditions, 
where 
$\boldsymbol{f}:Q \to \mathbb{R}^3$, 
$h:Q \to \mathbb{R}^{3\times3}_{\textrm{sym}}$, and 
${\risei \sigma_{*}}:Q\to \mathbb{R}^{3\times3}_{\textrm{sym}}$ are given functions in $Q$, 
{\fukao $\mathbb{R}^{3\times3}_{\textrm{sym}}$ stands for the $3 \times 3$ symmetric matrix. }
With the help of $h$ and $\sigma _*$, 
we can translate the problem to the homogeneous boundary value problem.  
The operator $\mathbf{div}$ is defined by 
$\mathbf{div} \tau:=(\textrm{div}\tau_{1\cdot},\textrm{div}\tau_{2\cdot},\textrm{div}\tau_{3\cdot})$ 
{\fukao for all $\tau \in \mathbb{R}^{3\times3}_{\textrm{sym}}$, where 
$\textrm{div}\tau_{i\cdot}:=\sum _{j=1}^3 \partial \tau _{ij}/\partial x_j$ for 
$i=1,2,3$.}
The first equation \eqref{dleqv} is derived by the conservation law 
of momentum. 
The second equation (\ref{dleqs}) ensures the property of perfect plasticity, 
where we assume 
the sum decomposition of strain. 
In relation \eqref{dleqs}, 
${ Z}$ is a time-dependent closed convex set defined by 
${ Z} := \{\tau \in \mathbb{R}^{3\times3}_{\textrm{sym}} : (1/2)|\tau ^{D}|^2 \le { g} \}$, 
where 
$\tau^D_{ij}:=
\tau_{ij}-(1/3)\sum_{k=1}^3\tau_{kk}
\delta_{ij}$ for $i,j=1,2,3$ and 
$|\tau|^2:=
\sum_{i,j=1}^3 \tau_{ij}\tau_{ij}$ for all {\fukao $\tau \in \mathbb{R}^{3\times3}_{\textrm{sym}}$, 
the symbol $\delta_{ij}$ is the {K}ronecker delta.}
Moreover, $g$ is a positive constant that represents the threshold value of the 
elasto-plasticity. 
$\partial I_{ Z}$ is the subdifferential of $I_{ Z}$. 
{\fukao Originally, one considered 
\begin{align*}
	&\frac{\partial \hat{\boldsymbol{v}} }{\partial t} 
	= \mathbf{div} \hat{\sigma} + \hat{\boldsymbol{f}} \quad \mbox{in}~Q,\\
	&\frac{\partial \hat{\sigma}}{\partial t}+\partial I_{ Z}(\hat{\sigma})\ni \varepsilon 
	 (\hat{\boldsymbol{v}}) \quad \mbox{in}~Q, \\
	& \hat{\boldsymbol{v}}=\boldsymbol{h}_D \quad {\rm on }(0,T) \times \Gamma _D, \\
	& \hat{\sigma} \boldsymbol{\nu }=\boldsymbol{h}_N \quad {\rm on~}(0,T) \times \Gamma _N
\end{align*}
with some given enough smooth data $\hat{\boldsymbol{f}}$, $\boldsymbol{h}_D$ and
$\boldsymbol{h}_N$. 
Thanks to some suitable functions $\boldsymbol{f}$ and $h$, these can be translated to 
homogeneous boundary conditions. 
Indeed, follows from \cite[p.\ 238]{DL76} we use new variables 
$\boldsymbol{v}:=\hat{\boldsymbol{v}}-\boldsymbol{v}_*$ and 
$\sigma:=\hat{\sigma }-\sigma_*$ 
where $\boldsymbol{v}_*=\boldsymbol{h}_D$ 
on $(0,T) \times \Gamma _D$ and 
$\sigma_*\boldsymbol{\nu }=\boldsymbol{h}_N$ on 
$(0,T) \times \Gamma _N$, respectively. Then, we take 
$\boldsymbol{f}:=\hat{\boldsymbol{f}}-\partial \boldsymbol{v}_*/\partial t+\mathbf{div} \sigma_* $, 
and $h:=\varepsilon (\boldsymbol{v}_*)-\partial \sigma_* /\partial t$. 
The assumptions for $\boldsymbol{h}_D$ and $\boldsymbol{h}_N$ will be replaced by 
one for $\sigma _*$ and $h$ later.}\\

In this paper, we consider the problem in which the
threshold constant $g$ is replaced by the function $g:[0,T]\times\Omega\to(0,\infty)$.
{\fukao In the forthcoming paper, we will} focus on the case of an unknown dependent threshold $g(\sigma )$, 
which is a more 
interesting setting. Then, \eqref{dleqs} becomes the quasivariational inequality (see, e.g.\ 
\cite{BBW15, CR06, Ste04}). {\fukao See also interesting related topics \cite{Kis17, KK16}. }

\subsection{Notation}

Hereafter, we use the following notation: 
$\boldsymbol{H}:=L^2(\Omega )^3$, 
$\boldsymbol{V}:=\{ \boldsymbol{z} \in H^1(\Omega )^3 : \boldsymbol{z}=\boldsymbol{0} \ 
\mbox{a.e.\ } \mbox{on\ } \Gamma _D\}$, 
with their inner products $(\cdot, \cdot )_{\boldsymbol{H}}$, 
$(\cdot, \cdot )_{\boldsymbol{V}}$, and the norm $| \cdot |_{\boldsymbol{H}}$, where 
$| \cdot |_{\boldsymbol{V}}$ is defined by
\begin{equation*} 
	|\boldsymbol{z}|_{\boldsymbol{V}}
	:= \left\{ 
	\sum_{i,j=1}^3 \int_{\Omega} 
	\left|\frac{\partial z_i}{\partial x_j} \right|^2 dx\right\}^{\frac{1}{2}}
	\quad \mbox{for~all~} \boldsymbol{z} \in \boldsymbol{V}. 
\end{equation*} 
Denote the dual space of $\boldsymbol{V}$ by $\boldsymbol{V}^*$ with 
the duality pair $\langle \cdot ,\cdot \rangle _{\boldsymbol{V}^*,\boldsymbol{V}}$. 
Moreover, we define the following bilinear form: 
$(\!( \cdot, \cdot )\!) :\boldsymbol{V} \times \boldsymbol{V} \to \mathbb{R}$
\begin{equation*} 
	(\!( \boldsymbol{z}, \tilde{\boldsymbol{z}} )\!)
	:=\sum_{i,j=1}^3\int_{\Omega} \frac{\partial z_i}{\partial x_j} \frac{\partial \tilde{z}_i}{\partial x_j} dx
	\quad \mbox{for~all~} \boldsymbol{z} ,\tilde{\boldsymbol{z}} \in \boldsymbol{V}.
\end{equation*}
We also define $\mathbb{H}:=\{ \tau:=\{\tau_{ij}\} : \tau_{ij}\in L^2(\Omega),\ \tau_{ij}=\tau_{ji}\}$, 
$\mathbb{V}:=\{ \tau \in \mathbb{H} : \textbf{div} \tau \in \boldsymbol{H}, \tau _{i \cdot } \cdot \boldsymbol{\nu } =0 
\ \mbox{a.e.\ on\ } \Gamma _N \}$ with their inner products 
\begin{align*} 
	(\tau ,\tilde{\tau} )_{\mathbb{H}} 
	& :=\sum_{i,j=1}^3\int_{\Omega}\tau_{ij} \tilde{\tau}_{ij}dx
	\quad \mbox{for~all~}\tau ,\tilde{\tau } \in \mathbb{H},
	\\ 
	(\tau ,\tilde{\tau} )_{\mathbb{V}}
	& := (\tau ,\tilde{\tau} )_{\mathbb{H}} 
	+(\mathbf{div}\tau ,\mathbf{div}\tilde{\tau } )_{\boldsymbol{H}}
	=(\tau ,\tilde{\tau} )_{\mathbb{H}} 
	+\sum_{i,j=1}^3 \int_{\Omega} 
	\frac{\partial \tau_{ij}}{\partial x_j}
	\frac{\partial \tilde{\tau}_{ij}}{\partial x_j}
	 dx 
	 \quad \mbox{for~all~}\tau ,\tilde{\tau } \in \mathbb{V}.
\end{align*}
The following convex constraint plays an important role in this paper. 
For each $t\in [0,T]$, 
\begin{equation*} 
	\tilde{K}(t):=
	\left\{ 
	\tau \in \mathbb{H} : 
	\frac{1}{2} \bigl| \tau^D(x) \bigr|^2 \le g(t,x) \quad \mbox{for~a.a.\ } x \in \Omega \right\}
	{\fukao ,}
	\quad K(t):= \tilde{K}(t)- \sigma_*(t).
\end{equation*}

Finally, we recall an important relation. 
For each $\boldsymbol{z} \in \boldsymbol{V}$, $\tau \in \mathbb{V}$, 
the following relation holds: 
\begin{equation}
	\bigl( 
	\varepsilon (\boldsymbol{z}),
	\tau
	\bigr)_{\mathbb{H}}
	+
	(\mathbf{div} \tau, \boldsymbol{z}
	)_{\boldsymbol{H}}=0.
	\label{gg}
\end{equation}
This is called the {G}auss--{G}reen relation.

\subsection{Definition of the solution and theorem}

Under these settings, 
we {\fukao define a solution of modified problem for \eqref{dleqv} and \eqref{dleqs}} in the variational formulation:

\paragraph{Definition 2.1.} 
\textit{For each {\fukao $\kappa \in (0,1]$ and $\nu \in (0,1]$}, the pair $(\boldsymbol{v}, \sigma)$ is called a solution 
{\fukao of modified problem for \eqref{dleqv} and \eqref{dleqs}}
in the sense of variational inequality if
\begin{gather*}
	\boldsymbol{v} \in H^{1}(0,T;\boldsymbol{H}) \cap L^\infty (0,T;\boldsymbol{V}) \cap 
	L^2 \bigl( 0,T;H^2(\Omega )^3 \bigr),
	\\
	\sigma \in H^1(0,T;\mathbb{H}) \cap L^2(0,T;\mathbb{V}), 
	\quad \sigma (t) \in K(t) \quad \textit{ for~all~} t \in [0,T],
\end{gather*}
and they satisfy
\begin{gather*}
	\bigl(
	\boldsymbol{v}'(t),
	\boldsymbol{z}
	\bigr)_{\boldsymbol{H}}
	+ 
	\nu \bigl( \! \bigl( 
	\boldsymbol{v}(t), 
	\boldsymbol{z} 
	\bigr) \! \bigr) 
	- \bigl( 
	\mathbf{div} \bigl( \sigma(t) \bigr), \boldsymbol{z}
	\bigr)_{\boldsymbol{H}}
	= 
	\bigl( 
	\boldsymbol{f}(t),\boldsymbol{z} 
	\bigr)_{\boldsymbol{H}} 
	\quad 
	\textit{for~all~}
	\boldsymbol{z} \in \boldsymbol{V}, \label{vi1} \\
	\bigl( 
	\sigma'(t), \sigma(t)- \tau 
	\bigr)_{\mathbb{H}} 
	+ \kappa 
	\bigl( \sigma (t), \sigma (t)-\tau \bigr)_{\mathbb{V}}
	-
	\bigl( 
	\varepsilon \bigl( 
	\boldsymbol{v}(t)
	\bigr), 
	\sigma(t)-\tau 
	\bigr)_{\mathbb{H}}
	\notag \\
	\le 
	\bigl( 
	h(t), 
	\sigma(t)-\tau 
	\bigr)_{\mathbb{H}}
	\quad 
	\textit{for~all~}
	\tau \in K(t)\cap \mathbb{V} \label{vi2} 
\end{gather*}
for a.a.\ $t \in (0,T)$ with $\boldsymbol{v}(0)=\boldsymbol{v}_0$ 
in $\boldsymbol{H}$ and $\sigma(0)=\sigma_0$ in $\mathbb{H}$.} \\

{\fukao As the remark, $\kappa \in (0,1]$ make no physical meanings, therefore we should consider the case that $\kappa =0$ (see, {D}efinition~2.3). 
On the other hand, $\nu \in (0,1]$ means the viscosity coefficient for some damping. } 
 
Hereafter, we assume that 
\begin{enumerate}
 \item[(A1)] $\boldsymbol{f} \in L^2(0,T;\boldsymbol{H})$ and $h \in L^2(0,T;\mathbb{H})$;
 \item[(A2)] $\boldsymbol{v}_0 \in \boldsymbol{V}$ and $\sigma _0 \in K(0) \cap\mathbb{V}$;
 \item[(A3)] $\sigma _* \in H^1(0,T;\mathbb{V})$;
 \item[(A4)] {\fukao $g \in H^1(0,T;C( \overline{\Omega})) \cap C(\overline{Q})$};
 \item[(A5)] There exist two constants $C_1, C_2>0$ such that 
\begin{equation*} 
	0 < C_1 \le g(t,x) \le C_2 \quad \mbox{for~all~} (t,x) \in {\fukao \overline{Q}}.
\end{equation*} 
\end{enumerate}

From the definition of $K(t)$, we see that 
$-\sigma _* (t) \in K(t)$ for all $t \in [0,T]$. 
Our first theorem is as follows:

\paragraph{Theorem 2.1.} 
\textit{Under assumptions {\upshape (A1)--(A5)}, there exists a unique solution 
$(\boldsymbol{v}, \sigma)$ 
{\fukao of modified problem for \eqref{dleqv} and \eqref{dleqs}}
in the sense of variational inequality.} \\

To relax assumption (A4) on $g$ with respect to time regularity, 
we recall the concept of the weak variational formulation:

\paragraph{Definition 2.2.} 
\textit{For each {\fukao $\kappa \in (0,1]$ and $\nu \in (0,1]$}, 
the pair $(\boldsymbol{v}, \sigma)$ is called a solution 
{\fukao of modified problem for \eqref{dleqv} and \eqref{dleqs}}
in the sense of weak variational inequality if
\begin{gather*}
	\boldsymbol{v} \in H^{1}(0,T;\boldsymbol{H}) \cap L^\infty (0,T;\boldsymbol{V}) \cap 
	L^2 \bigl( 0,T;H^2(\Omega )^3 \bigr),
	\\
	\sigma \in C\bigl([0,T];\mathbb{H} \bigr) \cap L^2(0,T;\mathbb{V}), 
	\quad \sigma (t) \in K(t) \quad \textit{ for~a.a.\ } t \in [0,T],
\end{gather*}
and they satisfy
\begin{gather}
	\bigl(
	\boldsymbol{v}'(t),
	\boldsymbol{z}
	\bigr)_{\mathbb{H}}
	+ 
	\nu \bigl( \! \bigl( 
	\boldsymbol{v}(t), 
	\boldsymbol{z} 
	\bigr) \! \bigr) 
	- \bigl( \mathbf{div}
	\sigma(t), \boldsymbol{z}
	\bigr)_{\boldsymbol{H}}
	= 
	\bigl( 
	\boldsymbol{f}(t),\boldsymbol{z} 
	\bigr)_{\boldsymbol{H}} 
	\quad 
	\textit{for~all~}
	\boldsymbol{z} \in \boldsymbol{V}, \label{wvi1} \\
	\int_{0}^{t} \bigl( 
	\eta'(s), \sigma(s)- \eta(s) 
	\bigr)_{\mathbb{H}} ds 
	+\kappa \int_{0}^{t} \bigl( 
	\sigma(s), \sigma(s)- \eta(s) 
	\bigr)_{\mathbb{V}} ds 
	\notag \\
	 - 
	\int_{0}^{t}
	\bigl( 
	\varepsilon \bigl( 
	\boldsymbol{v}(s)
	\bigr), 
	\sigma(s)-\eta(s)
	\bigr)_{\mathbb{H}} ds 
	+ \frac{1}{2} \bigl| \sigma (t) - \eta (t) \bigr|_{\mathbb{H}}^2 
	\notag \\
	{}
	\le \int_{0}^{t}
	\bigl( 
	h(s), 
	\sigma(s)-\eta(s) 
	\bigr)_{\mathbb{H}} ds 
	+ \frac{1}{2}  \bigl| \sigma_0 - \eta (0) \bigr|_{\mathbb{H}}^2 
	\quad 
	\textit{for~all~}
	\eta \in {\mathcal K}_0,
	\label{wvi2}
\end{gather}
for a.a.\ $t \in (0,T)$ with $\boldsymbol{v}(0)=\boldsymbol{v}_0$ 
in $\boldsymbol{H}$ and $\sigma(0)=\sigma_0$ in $\mathbb{H}$,
where
\begin{equation*} 
	{\mathcal K}_0:=\bigl\{ 
	\eta \in H^1(0,T;\mathbb{H}) \cap L^2(0,T;\mathbb{V}):
	\eta (t) \in K(t) \quad \textit{for~a.a.\ } t \in (0,T) \bigr\}.
\end{equation*} 
}

We assume the weaker condition {\fukao (A4$'$)} in place of (A4):
\begin{enumerate}
 \item[{\fukao (A4$'$)}] $g \in C(\overline{Q})$.
\end{enumerate}

\paragraph{Theorem 2.2.} 
\textit{Under assumptions {\upshape (A1)--(A3)}, {\fukao {\upshape (A4$'$)}}, and {\upshape (A5)}, 
there exists a unique solution 
$(\boldsymbol{v}, \sigma)$ 
{\fukao of modified problem for \eqref{dleqv} and \eqref{dleqs}}
in the sense of weak variational inequality.} \\

Let us neglect the parameter {\fukao $\kappa\in (0,1]$}. In this case, the problem is the same as the \emph{{M}oreau sweeping process}.

\paragraph{Definition 2.3.} 
\textit{For each {\fukao $\nu \in (0,1]$}, 
the pair $(\boldsymbol{v}, \sigma)$ is called a solution of the {\fukao viscous}
perfect plasticity model for 
{\fukao \eqref{dleqv} and \eqref{dleqs}}
if
\begin{gather*}
	\boldsymbol{v} \in H^{1}(0,T;\boldsymbol{V}^*) \cap L^\infty (0,T;\boldsymbol{H}) \cap 
	L^2 \bigl( 0,T;\boldsymbol{V} \bigr),
	\\
	\sigma \in H^1(0,T;\mathbb{H}), 
	\quad \sigma (t) \in K(t) \quad \textit{ for~all~} t \in [0,T],
\end{gather*}
and they satisfy
\begin{gather*}
	\bigl\langle  
	\boldsymbol{v}'(t),
	\boldsymbol{z}
	\bigr\rangle _{\boldsymbol{V}^*, \boldsymbol{V}} 
	+ 
	\nu \bigl( \! \bigl( 
	\boldsymbol{v}(t), 
	\boldsymbol{z} 
	\bigr) \! \bigr) 
	+ \bigl( 
	\sigma(t),\varepsilon (\boldsymbol{z})
	\bigr)_{\mathbb{H}}
	= 
	\bigl( 
	\boldsymbol{f}(t),\boldsymbol{z} 
	\bigr)_{\boldsymbol{H}} 
	\quad 
	\textit{for~all~}
	\boldsymbol{z} \in \boldsymbol{V}, \label{pp1} \\
	\bigl( 
	\sigma'(t), \sigma(t)- \tau 
	\bigr)_{\mathbb{H}} 
	-
	\bigl( 
	\varepsilon \bigl( 
	\boldsymbol{v}(t)
	\bigr), 
	\sigma(t)-\tau 
	\bigr)_{\mathbb{H}}
	\le 
	\bigl( 
	h(t), 
	\sigma(t)-\tau 
	\bigr)_{\mathbb{H}}
	\quad 
	\textit{for~all~}
	\tau \in K(t) \label{pp2} 
\end{gather*}
for a.a.\ $t \in (0,T)$ with $\boldsymbol{v}(0)=\boldsymbol{v}_0$ 
in $\boldsymbol{H}$ and $\sigma(0)=\sigma_0$ in $\mathbb{H}$.} \\

We replace (A3) by {\fukao (A3$'$)}:
\begin{enumerate}
 \item[{\fukao (A3$'$)}] $\sigma _* \in H^1(0,T;\mathbb{H})$.
\end{enumerate}

Our last theorem is as follows:

\paragraph{Theorem 2.3.} 
\textit{Under assumptions {\upshape (A1), (A2), {\fukao (A3$'$)}, (A4)}, and {\upshape (A5)}, 
there exists a unique solution 
$(\boldsymbol{v}, \sigma)$ of the 
{\fukao viscous}
perfect plasticity model for 
{\fukao \eqref{dleqv} and \eqref{dleqs}}.}

\section{Proof of the main theorems}
\setcounter{equation}{0}

In this section, we recall the known results with respect to the 
time-dependent evolution inclusion. 
Moreover, we apply them to prove our main theorems.

\subsection{Evolution inclusion with time-dependent domain}

In this subsection, we recall the result given by {K}enmochi \cite{Ken81} and 
apply it to the proof of Theorem~2.1. 
For this purpose, we describe the solvability of the evolution inclusion 
generated by the subdifferential operator with a time-dependent domain. 
For the family $\{ \phi^t \}:=\{\phi^t\}_{t \in [0,T]}$ of time-dependent, proper,
lower semicontinuous, and convex functionals 
on {\risei the real Hilbert space} $H$ {\fukao equipped with the norm $|\cdot|_H$}, let us consider the following form: 
\begin{equation}
	u'(t)+\partial \phi^t \bigl( u(t) \bigr) \ni f(t)
	\quad \mbox{in~} H, \ \mbox{for~a.a.\ } t \in (0,T),
	\quad u(0)=u_0 \quad \mbox{in } ~H, 
	\label{kenmochi}
\end{equation}
where $f \in L^2(0,T;H)$ and $u_0 \in H$ are given functions and 
$\partial \phi^t$ is the subdifferential of $\phi^t$ in $H$. 
We introduce a proposition for the existence of solutions of \eqref{kenmochi} 
under the following condition (H) for $\{ \phi^t \}$:

\begin{quote}
(H)\quad  For each $r>0$, there exist 
$\alpha_r\in L^2(0,T)$ and $\beta_r\in L^1(0,T)$ 
satisfying the following property: For each $s, t \in [0,T]$ with $s \le t$ and 
$z\in D(\phi^s)$ with $| z |_H \le r$, there exists $\tilde{z}\in D(\phi^t)$ such that  
\begin{equation}
	|\tilde{z}-z|_H \le \left( \int_s^t \alpha_r(l) dl \right) \left( 
	1+{\fukao \bigl| \phi^s (z) \bigr|^{\frac{1}{2}} } \right) 
	\label{H1}
\end{equation}
and
\begin{equation}
	\phi^t (\tilde{z})-\phi^s(z) \le \left( \int_s^t\beta_r(l)dl\right) 
	\bigl( 1+ {\fukao \bigl| \phi^s(z) \bigr| } \bigr).
	\label{H2}
\end{equation}
\end{quote}

\paragraph{Proposition 3.1. \cite[{K}enmochi]{Ken81}} 
\textit{Assume that $\{ \phi^t \}$ satisfies condition
{\upshape (H)}. 
Then, for each $u_0 \in \overline{D(\phi^0)}$, {\fukao the closure of $D(\phi^0)$ with respect to $H$-norm},
and $f\in L^2(0,T;H)$, there exists a unique $u\in C([0,T];H)$ with
\begin{equation*} 
	\sqrt{t} u' \in L^2(0,T;H), \quad \sup_{t \in [0,T]} t \phi^t \bigl( u(t) \bigr) <\infty.
\end{equation*}
such that $u$ satisfies \eqref{kenmochi}. 
Moreover, if $u_0\in D(\phi^0)$, then
\begin{equation*}
	u'\in L^2(0,T;H),\quad \sup_{t \in [0,T]} \phi^t \bigl( u(t) \bigr)<\infty.
\end{equation*}
}

To apply the above proposition to our problem, 
we choose the functionals $\psi :\boldsymbol{H} \to [0, \infty ]$ and
$\varphi :\mathbb{H} \to [0,\infty ]$ given by
\begin{equation*}
	\psi (\boldsymbol{z}) := 
	\begin{cases}
	\displaystyle \frac{\nu}{2} 
	(\!( \boldsymbol{z},\boldsymbol{z} )\!) 
	 & \mbox{if~} \boldsymbol{z} \in \boldsymbol{V}, \\
	\infty & \mbox{if~} \boldsymbol{z} \in \boldsymbol{H} \setminus \boldsymbol{V},
	\end{cases} 
	\quad 
	\varphi (\tau) := 
	\begin{cases}
	\displaystyle \frac{\kappa}{2} 
	|\tau|_{\mathbb{V}}^2
	 & \mbox{if~} \tau \in \mathbb{V}, \\
	\infty & \mbox{if~} \tau \in {\risei \mathbb{H}\setminus \mathbb{V} }. 
	\end{cases} 
\end{equation*}
Moreover, for each $t \in [0,T]$, we define $\varphi ^t:=\varphi +I_{K(t)}$, where $I_{K(t)}$ is the 
indicator function on $K(t)$.
Now, it is possible to reconstruct 
our problem as the following system of {C}auchy problems.
\begin{gather}
	\boldsymbol{v}'(t) 
	+ 
	\partial \psi \bigl( 
	\boldsymbol{v}(t) 
	\bigr) 
	+ 
	E_1 \sigma(t) 
	= \boldsymbol{f}(t)
	\quad \mbox{in~} 
	{\fukao \boldsymbol{H}}, 
	\ \mbox{for~a.a.\ } t \in (0,T), 
	\label{CP1}
	\\
	\sigma'(t) 
	+ \partial \varphi ^t \bigl( \sigma(t) \bigr) 
	+ E_2 \boldsymbol{v}(t) 
	\ni 
	h(t)\quad \mbox{in~} 
	\mathbb{H}, 
	\ \mbox{for~a.a.\ } t \in (0,T), 
	\label{CP2}
	\\
	\boldsymbol{v}(0)
	=\boldsymbol{v}_0 \quad 
	\mbox{in~} \boldsymbol{H},
	\quad 
	\sigma(0)=\sigma_0 \quad 
	\mbox{in~} \mathbb{H},\label{CP3}
\end{gather}
where  
it is obvious that 
$D(\partial \psi )=\boldsymbol{V} \cap H^2(\Omega )^3$ and 
$\partial \psi (\boldsymbol{z})=-\Delta \boldsymbol{z}$. 
Moreover, 
$E_1: \mathbb{V} \to \boldsymbol{H}$ 
is defined by $E_1 \tau:= -\mathbf{div} \tau$
for all $\tau \in \mathbb{V}$, and 
$E_2: \boldsymbol{V} \to \mathbb{H}$ is 
$E_2 \boldsymbol{z}:=-\varepsilon (\boldsymbol{z})$ for all $\boldsymbol{z} \in \boldsymbol{V}$.

We have the following lemma:

\paragraph{Lemma 3.1.}
\textit{For each $t\in(0,T)$, 
the function
$\varphi^t$ is proper, lower semicontinuous, 
and convex on $\mathbb{H}$, and 
$\partial \varphi ^t$ is characterized as follows: 
For each $\tau  \in D(\varphi ^t)=K(t) \cap \mathbb{V}$, 
$\tau ^* \in \partial \varphi ^t(\tau )$ in $\mathbb{H}$ if and only if 
\begin{equation}
	(\tau ^*, \tilde{\tau }-\tau )_{\mathbb{H}} 
	\le \kappa (\tau ,\tilde{\tau }-\tau )_{\mathbb{V}}
	\quad \mbox{for~all~} \tilde{\tau } \in K(t)\cap \mathbb{V}.
	\label{chara}
\end{equation}
Moreover, {\fukao there exists $T_* \in (0,T]$ such that}
$\{ \varphi ^t\}$ satisfies condition {\upshape (H)} {\fukao for all $s, t \in (0,T)$ with $|t-s|<T_*$}. }

\paragraph{Proof.}
First, we show that $\varphi^t$ 
is a lower semicontinuous functional on $\mathbb{H}$.
Let $\{ \tau_n \}_{n \in \mathbb{N}} 
\subset D(\varphi ^t)=K(t)\cap \mathbb{V}$; 
$\tau_n \to \tau$ in $\mathbb{H}$ as $n \to \infty$. 
If $R:=\liminf_{n\to\infty} \varphi^t(\tau_n)=\infty$, 
then it is obvious. 
Thus, we assume $R<\infty$. 
We can take a subsequence $\{ \tau_{n_k} \}_{k \in \mathbb{N}}$ 
of $\{ \tau_n \}_{n \in \mathbb{N}}$ such that 
\begin{equation*}
	\varphi^t(\tau_{n_k}) \to R \quad \mbox{in }\mathbb{R}, \quad \tau_{n_k} \to \tau \quad \mbox{weakly in }
	\mathbb{V}, \ \mbox{as } k \to \infty.
\end{equation*}
As $K(t)\cap \mathbb{V}$ is a convex and 
closed set in $\mathbb{V}$, {\fukao that is, weakly closed from the convexity}.   
$\tau$ is surely an element of the set $K(t)\cap \mathbb{V}$. 
Therefore, by the lower semicontinuity of {\risei $\varphi$},    
\begin{align*}
	\liminf_{n\to\infty} \varphi^t(\tau_n) = {} & 
	\liminf_{k\to\infty} \varphi^t(\tau_{n_k}) \\
	= {} & \liminf_{k\to\infty} \bigl( {\fukao \varphi(\tau_{n_k})} + I_{K(t)} (\tau_{n_k}) 
	\bigr)\\
	\ge {} & {\fukao \varphi(\tau)}+I_{K(t)}(\tau)\\
	= {} & \varphi^t(\tau).
\end{align*}

Next, 
for each $\tau  \in D(\varphi ^t)=K(t) \cap \mathbb{V}$, 
take $\tau ^* \in \partial \varphi ^t(\tau )$ in $\mathbb{H}$. 
From the definition of the subdifferential, 
\begin{equation*}
	(\tau ^*, \bar{\tau }-\tau )_{\mathbb{H}} 
	\le \frac{\kappa }{2} |\bar{\tau }|_{\mathbb{V}}^2 - 
	\frac{\kappa }{2} |\tau|_{\mathbb{V}}^2
	\quad \textrm{for~all~} \bar{\tau } \in K(t)\cap \mathbb{V}.
\end{equation*}
Here, for each $\tilde{\tau} \in K(t) \cap \mathbb{V}$ and $r \in (0,1)$, 
substituting $\bar{\tau }:=r \tilde{\tau} +(1-r)\tau \in K(t) \cap \mathbb{V}$, 
we have
\begin{equation*}
	r(\tau ^*, \tilde{\tau }-\tau )_{\mathbb{H}} 
	\le 
	\kappa r ( \tau, \tilde{\tau }-\tau)_{\mathbb{V}} 
	+
	\frac{\kappa }{2} r^2 ( \tilde{\tau}-\tau , \tilde{\tau}-\tau )_{\mathbb{V}}.
\end{equation*}
Dividing this by $r$ and letting $r \to 0$, we obtain \eqref{chara}. 

Finally, from the definition of $K(s)$, we see that  
for each $\tau \in K(s)\cap \mathbb{V}$, 
there exists 
$\tilde{\tau} \in \tilde{K}(s)$ such that 
$\tau=\tilde{\tau}-\sigma_{*}(s)$. 
Using assumptions (A4) and (A5), there exists ${\fukao T_*} \in (0,T]$ such that
\begin{equation*}
	\bigl| g(t)-g(s) \bigr|_{C(\overline{\Omega})}<C_1 
	\label{gasm1}
\end{equation*}
for all $s,t \in(0,T)$ with $|t-s|<{\fukao T_*}$.  
Now, we take 
\begin{equation*}
	\tilde{\tau}_{*}:=\left(1-\frac{1}{C_1} \bigl| g(t)-g(s) 
	\bigr|_{C(\overline{\Omega})} \right) \tilde{\tau}.
	\label{tautil}
\end{equation*}
Then, because 
$\tilde{\tau }=\tau +\sigma _*({\fukao s}) \in \mathbb{V}$ and $\tilde{\tau}\in \tilde{K}(s)$,
\begin{align*}
	\bigl| \tilde{\tau}_{*}^D 
	\bigr|^2 = {} 
	& \left(1-\frac{1}{C_1} \bigl| g(t)-g(s) \bigr|_{C(\overline{\Omega})} \right)^2
	\bigl| 
	\tilde{\tau}^D \bigr|^2\\
	\le {} 
	& 2 \left(1-\frac{1}{C_1} \bigl| g(t)-g(s) \bigr|_{C(\overline{\Omega})}\right) 
	g(s) \\
	\le {} 
	& 2g(s)-\frac{2g(s)}{C_1} \bigl| g(t)-g(s) \bigr|_{C(\overline{\Omega})} 
	\\
	\le {} 
	& 2 g(s)-2 g(t)+2g(t) - 2 \bigl| g(t)-g(s) \bigr|_{C(\overline{\Omega})} \\
	\le {} & 2g(t).
\end{align*}
Therefore, $\tilde{\tau }_{*} \in \tilde{K}(t) \cap \mathbb{V}$, that is, 
$\tau_{*}:=\tilde{\tau}_{*}-\sigma_{*}(t)$ is an element of $K(t)\cap \mathbb{V}$. 
Next, we show that $\{ \varphi ^t \}$ satisfies condition (H) {\fukao on some small time interval}. 
For each $\tau \in K(s)$, we take $\tau_{*}$ and $\tilde{\tau }_{*}$ 
to be the same as above, 
and obtain {\fukao 
\begin{align*}
	|\tau_{*}-\tau|_{\mathbb{H}} \le {} &
	|\tilde{\tau}_{*}-\tilde{\tau}|_{\mathbb{H}} 
	+ 
	\bigl| \sigma_{*}(t)-\sigma_{*}(s) \bigr|_{\mathbb{H}}
	\\
	\le {} & 
	\frac{1}{C_1}
	\bigl| 
	g(t)-g(s)
	\bigr|_{C(\overline{\Omega})} |\tilde{\tau }|_{\mathbb{H}}
	+ 
	\bigl| 
	\sigma_{*}(t)-\sigma_{*}(s)
	\bigr|_{\mathbb{H}} 
	\\
	\le {} &
	\int_s^t \left( \frac{1}{C_1}
	\bigl| g'(l) \bigr|_{C(\overline{\Omega})}{\fukao |\tilde{\tau }|_{\mathbb{H}}}
	+
	\bigl| 
	\sigma_{*}'(l) 
	\bigr|_{\mathbb{H}} 
	\right) dl
	\\
	\le {} &
	\int_s^t \left( \frac{1}{C_1}
	\bigl| g'(l) \bigr|_{C(\overline{\Omega})} \bigl| \tilde{\tau }-\sigma _*(s) \bigr|_{\mathbb{H}}
	+\frac{1}{C_1}
	\bigl| g'(l) \bigr|_{C(\overline{\Omega})} |\sigma _*(s)|_{\mathbb{H}}
	+
	\bigl| 
	\sigma_{*}'(l) 
	\bigr|_{\mathbb{H}} 
	\right) dl. 
\end{align*}
Therefore, we can 
take 
\begin{equation*}
	\alpha_r(\cdot ):=\alpha (\cdot ):=
	\frac{1}{C_1} \bigl| g'(\cdot ) \bigr|_{C(\overline{\Omega})}
	\left(\sqrt{\frac{2}{\kappa}} +|\sigma _*|_{C([0,T];\mathbb{H})} \right)
	+|\sigma_{*}'(\cdot )|_{\mathbb{H}}
\end{equation*}
for each $r>0$}, and then 
$\alpha \in L^2(0,T)$ and  
\begin{equation}
	|\tau_{*}-\tau|_{\mathbb{H}}
	\le {\fukao \left(  \int_{s}^{t} \alpha (l) dl \right) 
	\left(1+ \varphi ^s(\tau)^{\frac{1}{2}} \right)}, 
	\label{h1}
\end{equation}
that is, 
\eqref{H1} is satisfied {\fukao for all $s, t \in (0,T)$ with $|t-s|<T_*$}. 
Additionally, setting $\theta:=1-(1/C_1)|g(t)-g(s)|_{C(\overline{\Omega})}$, 
we have $\theta < 1$ and
\begin{gather*}
	\varphi^t(\tau_{*}) = 
	\frac{\kappa}{2}\theta^2 
	(\tilde{\tau},\tilde{\tau})_{\mathbb{V}}
	- \kappa \theta \bigl(\tilde{\tau},\sigma_{*}(t) \bigr)_{\mathbb{V}}
	+ \frac{\kappa}{2} \bigl(\sigma_{*}(t),\sigma_{*}(t) \bigr)_{\mathbb{V}},
	\\
	\varphi^s(\tau) = 
	\frac{\kappa}{2} 
	(\tilde{\tau},\tilde{\tau})_{\mathbb{V}}
	- \kappa \bigl(\tilde{\tau},\sigma_{*}(s) \bigr)_{\mathbb{V}}
	+ \frac{\kappa}{2} \bigl(\sigma_{*}(s),\sigma_{*}(s) \bigr)_{\mathbb{V}}
\end{gather*}
from the definition of $\varphi^t$ and $\varphi ^s$. 
Thus, we have the following inequality:
\begin{align*}
	& \varphi^t (\tau_{*})-\varphi^s(\tau) \\
	\le {} 
	& \kappa |\tilde{\tau}|_{\mathbb{V}}
	\bigl| \sigma_{*}(t)-\sigma_{*}(s) \bigr|_{\mathbb{V}}
	+\frac{\kappa}{C_1} \bigl| g(t)-g(s) \bigr|_{C(\overline{\Omega})}
	|\tilde{\tau}|_{\mathbb{V}} \bigl| \sigma_{*}(t) \bigr|_{\mathbb{V}}
	\\
	& {} + \frac{\kappa}{2}
	\left( \bigl| \sigma_{*}(t) \bigr|^2_{\mathbb{V}}
	-
	\bigl| \sigma_{*}(s) \bigr|^2_{\mathbb{V}} \right)
	\\
	\le {} & \kappa 
	|\tilde{\tau}|_{\mathbb{V}}
	\bigl| \sigma_{*}(t)-\sigma_{*}(s) \bigr|_{\mathbb{V}} 
	+ \frac{\kappa }{C_1} \bigl| g(t)-g(s) \bigr|_{C(\overline{\Omega})}
	|\tilde{\tau}|_{\mathbb{V}} 
	|\sigma_{*}|_{C([0,T];\mathbb{V})} 
	\\
	& {}
	+ \kappa |\sigma_{*}|_{C([0,T];\mathbb{V})} 
	\left( \bigl| \sigma_{*}(t) \bigr|_{\mathbb{V}} 
	- \bigl| \sigma_{*}(s) \bigr|_{\mathbb{V}} \right)\\
	\le {} & 
	\left(\bigl(1+|\sigma_{*}|_{C([0,T];\mathbb{V})} \bigr)\int_s^t \bigl| \sigma_{*}'(l) \bigr|_{\mathbb{V}}dl
	+\frac{1}{C_1} 
	|\sigma_{*}|_{C([0,T];\mathbb{V})}
	\int_s^t \bigl| g'(l) \bigr|_{C(\overline{\Omega})}~dl 
	\right) 
	\bigl( 1+ |\tilde{\tau }|_{\mathbb{V}} \bigr). 
\end{align*} 
Therefore, 
if we take 
\begin{equation*}
	{\fukao \beta _r(\cdot ):=\beta (\cdot ):=
	(1+|\sigma_{*}|_{C([0,T];\mathbb{V})})
	\bigl| \sigma_{*}'(\cdot) \bigr|_{\mathbb{V}}
	+
	\frac{2}{C_1} 
	|\sigma_{*}|_{C([0,T];\mathbb{V})}
	\bigl| g'(\cdot) \bigr|_{C(\overline{\Omega})},
	}
\end{equation*}
then $\beta _r \in L^1(0,T)$ and  
\eqref{H2} is satisfied {\fukao for all $s, t \in (0,T)$ with $|t-s|<T_*$}. \hfill $\Box$

\paragraph{Proposition 3.2.} 
\textit{Let $\kappa \in (0,1]$. 
For any given 
$\tilde{\boldsymbol{v}} 
\in L^2(0,T;\boldsymbol{V})$, 
there exists a unique solution 
$\sigma \in H^1(0,T;\mathbb{H}) \cap L^{{\fukao \infty }}(0,T;\mathbb{V})$, with 
$\sigma (t) \in K(t)$ for all $t \in [0,T]$, 
to the following problem.
\begin{gather}
	\sigma'(t) 
	+ \partial \varphi ^t
	\bigl( \sigma(t) \bigr) + 
	E_2 \tilde{\boldsymbol{v}}(t) 
	\ni h(t) 
	\quad \mbox{in~} \mathbb{H}, 
	\ \mbox{for~a.a.\ } t \in (0,T),
	\label{prosig1}\\
	\sigma(0)=
	\sigma_0 \quad \mbox{in~}
	\mathbb{H}.
	\label{prosig2}
\end{gather}
Moreover, there exist positive constants $M_1$ and $M_2$, 
independent of $\kappa$ and $\nu$, such that
\begin{gather} 
	\bigl| 
	\sigma(t)
	\bigr|^2_{\mathbb{H}} 
	+ \kappa \int_{0}^{t} \bigl| \sigma(s) \bigr|_{\mathbb{V}}^2 ds
	\le M_1 \left( 1+
	\int_{0}^{T}\bigl| \tilde{\boldsymbol{v}}(s) \bigr|_{\boldsymbol{V}}^2 ds \right),
	\label{prosig3}
	\\
	\int_{0}^{t} \bigl| \sigma'(s) \bigr|_{\mathbb{H}}^2 ds
	+ \kappa \bigl| \sigma(t) \bigr|_{\mathbb{V}}^2 
	\le  
	M_2  \left( 1+
	\int_{0}^{T}\bigl| \tilde{\boldsymbol{v}}(s) \bigr|_{\boldsymbol{V}}^2 ds \right)
	\label{prosig4}
\end{gather}
for all $t \in [0,T]$. }

\paragraph{Proof.}

As a result of Lemma~3.1, we can apply Proposition~3.1 to 
obtain {\fukao a local solution $\sigma \in H^1(0,T_*;\mathbb{H}) 
\cap L^\infty (0,T_*;\mathbb{V})$} of \eqref{prosig1}--\eqref{prosig2}. 

Next, we obtain the uniform estimate \eqref{prosig3}. 
By the characterization \eqref{chara} 
of $\partial \varphi^t$, \eqref{prosig1} is equivalent to the following 
inequality:
\begin{equation}
	\bigl( 
	\sigma'(s),\sigma(s)-\tau
	\bigr)_{\mathbb{H}}
	+
	\kappa 
	\bigl( 
	\sigma(s),\sigma(s)-\tau
	\bigr)_{\mathbb{V}}
	-
	\bigl( 
	\varepsilon
	\bigl(
	\tilde{\boldsymbol{v}}(s)
	\bigr),
	\sigma(s)-\tau 
	\bigr)_{\mathbb{H}} \le 
	\bigl(
	h(s),\sigma(s)-\tau 
	\bigr)_{\mathbb{H}}
	\label{sub}
\end{equation}
for all $\tau \in K(s)$ and for a.a.\ $s \in (0,{\fukao T_*})$. From 
assumption (A3), we can 
substitute $\tau:=-\sigma_{*}(s) \in K(s)$ in \eqref{sub}, and then 
we obtain
\begin{align}
	& \frac{1}{2} \frac{d}{ds} 
	\bigl| \sigma(s)+\sigma_{*}(s) \bigr|^2_{\mathbb{H}}
	+ \kappa 
	\bigl| \sigma(s) \bigr|^2_{\mathbb{V}}
	 \notag \\
	& \le 
	-\bigl( \sigma_{*}'(s),\sigma(s)+\sigma_{*}(s) \bigr)_{\mathbb{H}}
	+
	\bigl( \varepsilon \bigl( \tilde{\boldsymbol{v}} (s) \bigr), 
	\sigma(s)+\sigma_{*}(s) \bigr)_{\mathbb{H}} 
	+ \bigl( 
	h(s), \sigma(s)+\sigma_{*}(s)
	\bigr)_{\mathbb{H}}
	\notag \\
	& \quad {}
	-\kappa \bigl( \sigma(s),\sigma_{*}(s) \bigr)_{\mathbb{V}}
	 \notag \\
	& \le
	\frac{3}{2} \bigl| \sigma _{*}'(s) \bigr|_{\mathbb{H}}^2
	+
	\frac{3}{2} \bigl| \tilde{\boldsymbol{v}}(s) \bigr|_{\boldsymbol{V}}^2  
	+ \frac{3}{2} \bigl| h(s) \bigr|_{\mathbb{H}}^2 + \frac{1}{2} \bigl| 
	\sigma(s)+\sigma_{*}(s) \bigr|_{\mathbb{H}}^2
	+ \frac{\kappa }{2} \bigl| \sigma(s) \bigr|^2_{\mathbb{V}}
	+ \frac{\kappa }{2} \bigl|\sigma_{*}(s) \bigr|^2_{\mathbb{V}}
	, \label{uni1}
\end{align}  
{\fukao where we use}
\begin{align}
	{\fukao \bigl| \varepsilon (\boldsymbol{z}) \bigr|_{\mathbb{H}}^2 }
	& {\fukao {}= \frac{1}{4}\sum _{i,j=1}^3 \int_{\Omega }^{} 
	\left| \frac{\partial z_i}{\partial x_j}+
	\frac{\partial z_j}{\partial x_i} \right|^2 dx } 
	\notag \\
	& {\fukao {}\le  \frac{1}{2}\sum _{i,j=1}^3 \int_{\Omega }^{} 
	\left\{ \left| \frac{\partial z_i}{\partial x_j}\right|^2+\left| 
	\frac{\partial z_j}{\partial x_i} \right|^2 \right\} dx } 
	\notag 
	\\
	& {\fukao {}= |\boldsymbol{z}|_{\boldsymbol{V}}^2 
	\quad {\rm for~all~} \boldsymbol{z} \in \boldsymbol{V}. 
	\label{ddd}}
\end{align}
{\fukao Above \eqref{uni1} means that}
\begin{equation*}
	\frac{d}{ds} \bigl| 
	\sigma(s)+\sigma_{*}(s)
	\bigr|^2_{\mathbb{H}} \notag
	\le 
	3 \bigl| \sigma _{*}'(s) \bigr|_{\mathbb{H}}^2
	+ 3 \bigl| \tilde{\boldsymbol{v}}(s) \bigr|_{\boldsymbol{V}}^2  
	+ 3 \bigl| h(s) \bigr|_{\mathbb{H}}^2 
	+ \bigl|\sigma_{*}(s) \bigr|^2_{\mathbb{V}}
	+ \bigl| 
	\sigma(s)+\sigma_{*}(s)
	\bigr|^2_{\mathbb{H}}
\end{equation*}  
for a.a.\ $s \in (0,T)$ with $\kappa \in (0,1]$. 
Now, using the {G}ronwall inequality, we see that 
there exists a constant $\tilde{M}_1$ that depends on 
$|\sigma_0 |_{\mathbb{H}}$, $|\sigma_*(0) |_{\mathbb{H}}$,  
$|\sigma_*'|_{L^2(0,T;\mathbb{H})}$, 
$|h|_{L^2(0,T;\mathbb{H})}$, $|\sigma_*|_{L^2(0,T;\mathbb{V})}$, and $T$ such that
\begin{equation*} 
	\bigl| 
	\sigma(t)+\sigma_{*}(t)
	\bigr|^2_{\mathbb{H}} 
	\le \tilde{M}_1 \left( 1+
	\int_{0}^{T}\bigl| \tilde{\boldsymbol{v}}(s) \bigr|_{\boldsymbol{V}}^2 ds \right).
\end{equation*} 
for all {\fukao $t \in [0,T_*]$}. 
This implies that \eqref{prosig3} holds {\fukao for all $t \in [0,T_*]$} for some suitable constant $M_1>0$. 

Next, we recall the approximate problem of \eqref{prosig1} using the 
{M}oreau--{Y}osida regularization $\varphi^t_\lambda $ for $\lambda \in (0,1]$ as follows:
\begin{equation*} 
	\varphi^t_\lambda(\tau ):=
	\inf_{\tilde{\tau }\in \mathbb{H}} 
	\left\{ 
	\frac{1}{2\lambda }
	|\tau -\tilde{\tau } |_{\mathbb{H}}^2 
	+ \varphi ^t(\tilde{\tau })
	\right\} 
	= \frac{1}{2\lambda }|\tau -J_\lambda ^t\tau |_{\mathbb{H}}^2 
	+ \varphi ^t (J_\lambda ^t\tau),
\end{equation*} 
where $J_\lambda ^t$ is the resolvent 
$(I+\lambda \partial \varphi ^t)^{-1}$ of $\partial \varphi ^t$. 
{\fukao Then we know that $\partial \varphi _\lambda ^t(\tau )=(\tau -J_\lambda ^t \tau )/\lambda $.}
From Lemma~3.1 with \eqref{h1}, for each $s,t \in [0,T]$ {\fukao with $|t-s|< T_*$} and 
$\tau \in \mathbb{H}$, there exists $\tilde{\tau } \in K(t)\cap \mathbb{V}$ such that 
\begin{gather*}
	|\tilde{\tau }-J_\lambda ^s \tau |_{\mathbb{H}} 
	\le {\fukao \left(  \int_{s}^{t} \alpha (l) dl \right) 
	\left(1+ \varphi ^s(J_\lambda ^s \tau)^{\frac{1}{2}} \right)}, \\
	\varphi^t (\tilde{\tau})-\varphi^s(J_\lambda ^s \tau) 
	\le \left( \int_s^t \beta(l)dl\right) 
	\bigl( 1+\varphi^s(J_\lambda ^s\tau ) 
	\bigr).
\end{gather*}
Hence, 
\begin{align*}
	\varphi ^t_\lambda (\tau )
	-\varphi ^s_\lambda (\tau ) 
	\le {} & 
	\frac{1}{2\lambda } |\tau -\tilde{\tau }|_{\mathbb{H}}^2 
	+\varphi ^t (\tilde{\tau})
	- \frac{1}{2\lambda } |\tau -J_\lambda ^s \tau |_{\mathbb{H}}^2 
	- \varphi ^s (J_\lambda ^s \tau)
	\\
	= {} & \frac{1}{2\lambda } \bigl( |\tau -\tilde{\tau }|_{\mathbb{H}}
	+|\tau -J_\lambda ^s \tau |_{\mathbb{H}} \bigr) 
	\bigl( |\tau -\tilde{\tau }|_{\mathbb{H}}
	-|\tau -J_\lambda ^s \tau |_{\mathbb{H}} \bigr) 
	+\varphi ^t (\tilde{\tau})- \varphi ^s (J_\lambda ^s \tau)
	\\
	\le  {} & \frac{1}{2\lambda } 
	\bigl( |\tau -J_\lambda ^s \tau+J_\lambda ^s \tau-\tilde{\tau }|_{\mathbb{H}}
	+|\tau -J_\lambda ^s \tau |_{\mathbb{H}} \bigr) 
	|\tilde{\tau }-J_\lambda ^s \tau |_{\mathbb{H}} 
	\\
	& {} +\left( \int_s^t \beta(l)dl\right) 
	\bigl( 1+\varphi^s(J_\lambda ^s\tau ) \bigr)
	\\
	\le {} & |\tilde{\tau }-J_\lambda ^s \tau|_{\mathbb{H}} 
	\left| \frac{\tau -J_\lambda ^s \tau }{\lambda }\right|_{\mathbb{H}}
	+ \frac{1}{2\lambda }|\tilde{\tau } -J_\lambda ^s \tau |_{\mathbb{H}}^2 
	+
	\bigl( 1+\varphi^s(J_\lambda ^s\tau )  \bigr)
	\int_s^t \beta(l)dl
	\\
	\le {} & \bigl| \partial \varphi _\lambda ^s(\tau )\bigr|_{\mathbb{H}} 
	{\fukao \left(  \int_{s}^{t} \alpha (l) dl \right) 
	\left(1+ \varphi ^s_\lambda ( \tau)^{\frac{1}{2}} \right)}
	\\
	& {\fukao {} + \left( \frac{t-s}{2\lambda } \int_{s}^{t} \bigl| \alpha (l) \bigr|^2 dl \right)
	\times 2\bigl( 1+ \varphi ^s_\lambda ( \tau) \bigr)
	 +
	\left( \int_s^t \beta(l)dl \right) \bigl( 1+\varphi^s_\lambda (\tau )  \bigr),}
\end{align*}
for a.a.\ $s \in (0,T)$, because 
\begin{equation*} 
	\varphi^s(J_\lambda ^s\tau ) =
	- \frac{1}{2\lambda }|\tau -J_\lambda ^s\tau |_{\mathbb{H}}^2 
	+ \varphi _\lambda ^s (\tau)
	\le \varphi ^s_\lambda (\tau ).
\end{equation*} 
Then, we obtain 
\begin{equation*}
	\frac{d}{ds} \varphi _\lambda^s (\tau ) \le 
	\bigl| \partial \varphi _\lambda ^s(\tau )
	\bigr|_{\mathbb{H}} \bigl| \alpha (s) \bigr| {\fukao \left(1+ \varphi ^s_\lambda ( \tau)^{\frac{1}{2}} \right)}
	+  \bigl| \beta (s) \bigr| {\fukao \bigl( 1+\varphi^s_\lambda (\tau ) \bigr)}
\end{equation*}
for a.a.\ $s \in (0,T)$, namely 
\begin{align*}
	& \frac{d}{ds} \varphi _\lambda^s  ({\fukao \eta}(s)) 
	- \bigl({\fukao \eta}'(s), \partial \varphi _\lambda^s  \bigl( 
	{\fukao \eta}(s) \bigr) \bigr)_{\mathbb{H}} 
	\\
	& \le 
	{\fukao 
	\bigl| \partial \varphi _\lambda ^s(\eta )
	\bigr|_{\mathbb{H}} \bigl| \alpha (s) \bigr| \left(1+ \varphi ^s_\lambda ( \eta)^{\frac{1}{2}} \right) 
	+  \bigl| \beta (s) \bigr| \bigl( 1+\varphi^s_\lambda (\eta ) \bigr) }
\end{align*}
for a.a.\ $s \in (0,T)$ and for all 
${\fukao \eta} \in W^{1,1}(0,T;\mathbb{H})$
(see, e.g.\ \cite[Lemma~1.2.5.]{Ken81}). 
From this rigorous formulation, we can obtain an estimate by
multiplying {\fukao approximate level of} \eqref{prosig1} by 
$\partial \varphi^s_{\fukao \lambda } \bigl(\sigma _{{\fukao \lambda }}(s) \bigr)$ at time $t=s$. 
That gives the following estimate:
\begin{align}
	& \bigl| \partial \varphi^s_{\fukao \lambda } \bigl( \sigma_{{\fukao \lambda }}(s) \bigr) \bigr|_{\mathbb{H}}^2 
	+ \frac{d}{ds} \varphi^s_{\fukao \lambda } \bigl( \sigma_{{\fukao \lambda }}(s) \bigr) 
	\notag \\
	& \le  
	\bigl| \partial \varphi^s_{\fukao \lambda } \bigl(\sigma_{{\fukao \lambda }}(s) \bigr) \bigr|_{\mathbb{H}} 
	\bigl| \alpha (s) \bigr| {\fukao \left(1+ \varphi ^s_\lambda \bigl(\sigma_{{\fukao \lambda }}(s) \bigr)^{\frac{1}{2}} \right)} 
	+ \bigl| \beta (s) \bigr| {\fukao \bigl( 1+\varphi^s_\lambda  \bigl(\sigma_{{\fukao \lambda }}(s) \bigr) }
	\notag \\
	& \quad {}
	- \bigl( \varepsilon \bigl( \tilde{\boldsymbol{v}}(s) \bigr), \partial \varphi^s_{\fukao \lambda }
	\bigl( \sigma_{{\fukao \lambda }}(s) \bigr)
	\bigr)_{\mathbb{H}} 
	+
	\bigl( h(s),  \partial \varphi^s_{\fukao \lambda } \bigl( \sigma(s) \bigr)
	\bigr)_{\mathbb{H}} 
	\notag \\
	& \le  
	\frac{1}{2} \bigl| \partial \varphi^s_{\fukao \lambda }\bigl( \sigma_{{\fukao \lambda }}(s) \bigr) \bigr|_{\mathbb{H}}^2 +
	\frac{3}{2} \bigl| \alpha (s) \bigr|^2
	+ \frac{3}{2}
	\bigl| \tilde{\boldsymbol{v}}(s) \bigr|_{\boldsymbol{V}}^2 
	+ \frac{3}{2}
	\bigl|h(s) 
	\bigr|_{\mathbb{H}}^2
	+
	 \bigl| \beta (s) \bigr| {\fukao \bigl( 1+\varphi^s_\lambda  \bigl(\sigma_{{\fukao \lambda }}(s) \bigr) }
	\label{prim}
\end{align}
for a.a.\ $s \in (0,T)$. Actually, we need to consider the above estimate 
at the level of 
the approximate problem of \eqref{prosig1} using the {M}oreau--{Y}osida regularization. 
Using the {G}ronwall inequality, we get 
\begin{equation*}
	\varphi ^t_{\fukao \lambda }\bigl( \sigma_{{\fukao \lambda }}(t) \bigr) 
	\le  
	\left( \frac{1}{2} |\sigma _0|_{\mathbb{V}}^2 
	+
	\frac{3}{2} |\alpha|_{L^2(0,T)}^2
	+ 
	\frac{3}{2} |h|_{L^2(0,T;\mathbb{H})}^2 
	+ 
	|\beta |_{L^1(0,T)}
	+
	\frac{3}{2} \int_{0}^{T}
	\bigl| \tilde{\boldsymbol{v}}(s) \bigr|_{\boldsymbol{V}}^2 ds 
	\right) 
	e^{|\beta |_{L^1(0,T)}}
\end{equation*}
for all $t \in [0,T]$, with $\kappa \in (0,1]$. 
Integrating \eqref{prim} over $[0,t]$ with respect to time, and using the above resultant,
 we obtain
\begin{equation*}
	\int_{0}^{t} \bigl| \partial \varphi^s _{\fukao \lambda } \bigl( \sigma_{\fukao \lambda }(s)\bigr)
	 \bigr|_{\mathbb{H}}^2 ds
	+ \frac{\kappa}{2} \bigl| \sigma_{\fukao \lambda }(t) \bigr|_\mathbb{V}^2
	\le 
	\tilde{M}_2 \left( 1+
	\int_{0}^{T}\bigl| \tilde{\boldsymbol{v}}(s) \bigr|_{\boldsymbol{V}}^2 ds \right)
\end{equation*}
for all {\fukao $t \in [0,T_*]$}, 
where $\tilde{M}_2$ is 
a positive constant depending on 
${\risei |\sigma_0 |_{\mathbb{V}}}$, 
$|h|_{L^2(0,T;\mathbb{H})}$, $|\alpha |_{L^2(0,T)}$, $|\beta |_{L^1(0,T)}$, and $T$, 
but independent of $\kappa \in (0,1]$. 
By comparison, we obtain 
{\fukao 
\begin{gather*} 
	\int_{0}^{t} \bigl| \sigma'_\lambda (s) \bigr|_{\mathbb{H}}^2 ds
	+ \kappa \bigl| \sigma_\lambda (t) \bigr|_{\mathbb{V}}^2 
	\le  
	M_2  \left( 1+
	\int_{0}^{T}\bigl| \tilde{\boldsymbol{v}}(s) \bigr|_{\boldsymbol{V}}^2 ds \right)
\end{gather*}
for all $t \in [0,T_*]$ 
for some suitable constant $M_2>0$. Finally, executing the limiting procedure 
we obtain \eqref{prosig3}, \eqref{prosig4} and $\sigma (t) \in K(t)$ for all $t \in [0,T_*]$.
Moreover, we can choose $T_*$ independent of the initial data. 
Therefore, we can use same argument on $[T_*, 2T_*]$ with the initial 
data $\sigma (T_*) \in K(T_*) \cap \mathbb{V}$. By the finite iteration, we show the existence and uniform estimates \eqref{prosig3} and \eqref{prosig4} for all $t \in [0,T]$.  
}
\hfill $\Box$ \\

We now define a solution operator $S_1:L^2(0,T;\boldsymbol{V}) \to L^2(0,T;\mathbb{H})$. 
This assigns a unique solution $S_1 \tilde{\boldsymbol{v}}:=\sigma$ 
to the above problem \eqref{prosig1}--\eqref{prosig2} on $[0,T]$.

\paragraph{Proposition 3.3.}
\textit{Let $\nu \in (0,1]$. 
For any given $\tilde{\sigma} \in L^2(0,T;\mathbb{V})$, there exists a unique solution 
$\mbox{\boldmath $ v$} \in H^1(0,T;\boldsymbol{H}) 
\cap L^\infty (0,T;\boldsymbol{V}) \cap L^2(0,T;H^2(\Omega )^3)$ 
to the following problem:
\begin{gather}
	\boldsymbol{v}'(t) + \partial \psi \bigl( \boldsymbol{v}(t) \bigr) 
	+ E_1 \tilde{\sigma}(t) = \boldsymbol{f}(t) \quad \mbox{in~} \boldsymbol{H}, \ \mbox{for~a.a.~} t \in (0,T), 
	\label{prov1}
	\\
	\boldsymbol{v}(0)=\boldsymbol{v}_0 \quad \mbox{in~} \boldsymbol{H}.
	\label{prov2}
\end{gather}
Moreover, there exist positive constants $M_3$ and $M_4$ which are independent of $\nu $ such that
\begin{gather*} 
	\bigl| 
	\boldsymbol{v}(t)
	\bigr|^2_{\boldsymbol{H}} 
	+\nu \int_{0}^{t} \bigl| \boldsymbol{v}(s) \bigr|^2_{\boldsymbol{V}} ds
	\le M_3 \left( 1+
	\int_{0}^{T} \bigl| \tilde{\sigma }(s) \bigr|_{\mathbb{V}}^2 ds \right),
	\label{prov3}
	\\
	\int_{0}^{t}\bigl| 
	\boldsymbol{v}'(s)
	\bigr|^2_{\boldsymbol{H}} ds
	+\nu \bigl| \boldsymbol{v}(t) \bigr|^2_{\boldsymbol{V}} 
	\le M_4 \left( 1+
	\int_{0}^{T} \bigl| \tilde{\sigma }(s) \bigr|_{\mathbb{V}}^2 ds \right)
	\label{prov4}
\end{gather*}
for all $t \in [0,T]$. } \\

As the proof uses standard techniques {\fukao (see, e.g.\ \cite{Bre73})}, we omit it from this paper. \\

We can also define a solution operator $S_2:L^2(0,T;\mathbb{V})\to L^2(0,T;\boldsymbol{V})$ 
that assigns a unique solution $S_2 \tilde{\sigma}:=\boldsymbol{v}$ 
to the above problem \eqref{prov1}--\eqref{prov2} on $[0,T]$.

\subsection{Proof of Theorem 2.1.}

In this subsection, we prove Theorem~2.1.

\paragraph{Proof of Theorem~2.1.}

We define the operator $S:L^2(0,T;\boldsymbol{V}) \to L^2(0,T;\boldsymbol{V})$ 
by $S:=S_2\circ S_1$. 
Let $\kappa, \nu \in (0,1]$. 
We show that $S$ is a contraction operator. 
For $\tilde{\boldsymbol{v}}^{(i)} \in L^2(0,T;\boldsymbol{V})$ for $i=1,2$, 
we set $\sigma^{(i)}:=S_1 \tilde{\boldsymbol{v}}^{(i)}$. 
Taking the difference between \eqref{prosig1} with $\sigma^{(1)}$ and 
{\risei\eqref{prosig1}} with $\sigma^{(2)}$ at {\fukao time $t=s$}, and using 
\eqref{chara} and the {G}auss--{G}reen relation \eqref{gg}, we have 
\begin{align*}
	& \frac{1}{2} \frac{d}{ds} \bigl| \sigma^{(1)}(s)-\sigma^{(2)}(s) \bigr|_{\mathbb{H}}^2 
	+ \kappa \bigl| \sigma^{(1)}(s)-\sigma^{(2)}(s) \bigr|_{\mathbb{V}}^2 
	\\
	& \le - \bigl( \mathbf{div} \bigl(\sigma^{(1)}(s)-\sigma^{(2)}(s)\bigr), \tilde{\boldsymbol{v}}^{(1)}(s) - \tilde{\boldsymbol{v}}^{(2)}(s) \bigr) _{\fukao \boldsymbol{H}}
	\nonumber \\
	& \le 
	\frac{\kappa }{2} \bigl|  \sigma^{(1)}(s)-\sigma^{(2)}(s) \bigr|_{\mathbb{V}}^2
	+
	\frac{1}{2\kappa } \bigl|  \tilde{\boldsymbol{v}}^{(1)}(s) - \tilde{\boldsymbol{v}}^{(2)}(s) \bigr|_{\boldsymbol{H}}^2
\end{align*} 
for a.a.\ $s \in (0,T)$. Integrating this expression over $(0,t)$ with respect to time, we obtain
\begin{equation}
	\bigl| \sigma^{(1)}(t)-\sigma^{(2)}(t) \bigr|_{\mathbb{H}}^2 
	+ \kappa \int_{0}^{t} \bigl| \sigma^{(1)}(s)-\sigma^{(2)}(s) \bigr|_{\mathbb{V}}^2 ds
	\le \frac{1}{\kappa }
	\int_0^t
	\bigl| \tilde{\boldsymbol{v}}^{(1)}(s) 
	- \tilde{\boldsymbol{v}}^{(2)}(s) \bigr|_{\boldsymbol{H}}^2ds
	\label{banach2}
\end{equation}
for all $t \in [0,T]$. 

Next, we set $\boldsymbol{v}^{(i)}:=S_2 \sigma^{(i)}$ for 
$i=1,2$. 
In the same way, taking the difference between \eqref{prov1} with 
$\boldsymbol{v}^{(1)}$ and \eqref{prov1} with $\boldsymbol{v}^{(2)}$ at 
{\fukao time $t=s$}, and 
using the {G}auss--{G}reen relation \eqref{gg}, 
we have
\begin{align*}
	& \frac{1}{2} \frac{d}{ds} \bigl| \boldsymbol{v}^{(1)}(s)-\boldsymbol{v}^{(2)}(s) \bigr|^2_{\boldsymbol{H}}
	+ \nu \bigl| \boldsymbol{v}^{(1)}(s) - \boldsymbol{v}^{(2)}(s) \bigr|_{\boldsymbol{V}}^2 \notag \\
	& \le - \bigl( \varepsilon \bigl( \boldsymbol{v}^{(1)}(s)-\boldsymbol{v}^{(2)}(s) \bigr), \sigma ^{(1)}(s)-\sigma ^{(2)}(s) \bigr)_{\mathbb{H}} \notag \\
	& \le 
	\frac{\nu }{2} \bigl| \boldsymbol{v}^{(1)}(s)-\boldsymbol{v}^{(2)}(s) \bigr|_{\boldsymbol{V}}^2
	+
	\frac{1}{2\nu } \bigl| \sigma ^{(1)}(s)-\sigma ^{(2)}(s) \bigr|_{\mathbb{H}}^2 ,
\end{align*}
that is,
\begin{equation}
	\bigl| \boldsymbol{v}^{(1)}(t)-\boldsymbol{v}^{(2)}(t) \bigr|^2_{\boldsymbol{H}}
	+\nu \int_{0}^{t} \bigl| \boldsymbol{v}^{(1)}(s) - \boldsymbol{v}^{(2)}(s) \bigr|_{\boldsymbol{V}}^2 ds 
	\le \frac{1}{\nu } \int_{0}^{t} \bigl| \sigma ^{(1)}(s)-\sigma ^{(2)}(s) \bigr|_{\mathbb{H}}^2 ds
	\label{banach3}
\end{equation} 
for all $t \in [0,T]$. 
Combining \eqref{banach3} with \eqref{banach2}, we obtain
the following inequality:
\begin{align*}
	\nu \bigl| S \tilde{\boldsymbol{v}}^{(1)}
	-S \tilde{\boldsymbol{v}}^{(2)}\bigr|^2_{L^2(0,T;\boldsymbol{V})} 
	\le {} & \frac{1}{\kappa \nu }T  \bigl| 
	\tilde{\boldsymbol{v}}^{(1)} - \tilde{\boldsymbol{v}}^{(2)} 
	\bigr|_{L^2(0,T;\boldsymbol{V})}^2.
	\label{finalest}
\end{align*} 
Taking $T_0 \in (0,T]$ satisfying 
$(1/\kappa \nu^2)T_0 <1$, 
we see that $S$ is a contraction on $L^2(0,T_0;\boldsymbol{V})$.
Thus, we can apply the {B}anach fixed point theorem to show that 
there exists a unique  
$\boldsymbol{v} \in L^2(0,T_0;\boldsymbol{V})$ 
such that $\boldsymbol{v}=S \boldsymbol{v}$ in $L^2(0,T_0;\boldsymbol{V})$ at once. 
Then, the pair $(\boldsymbol{v}, \sigma )$ of $\boldsymbol{v}$ and 
$\sigma :=S_1 \boldsymbol{v}$ is 
the solution in the sense of variational inequality on $[0,T_0]$. 
Next, recall that we can choose $T_0$ 
independent of the initial data. 
Taking $\boldsymbol{v}(T_0) \in \boldsymbol{V}$ and 
$\sigma (T_0) \in K(T_0) \cap \mathbb{V}$ as the initial data at time $T_0$, 
we can find the solution pair $(\boldsymbol{v}, \sigma)$ in the sense of variational inequality on $[T_0,2T_0]$, 
i.e.\ we obtain the solution $(\boldsymbol{v}, \sigma)$ 
{\fukao of modified problem for \eqref{dleqv} and \eqref{dleqs}}
in the sense of variational inequality on $[0,T]$ by iterating at finite times. 
\hfill $\Box$

\subsection{Proof of Theorem 2.2.}

In this subsection, we prove Theorem~2.2. 
Assumption (A4) is replaced by {\fukao (A4$'$)}. 
From assumptions {\fukao (A4$'$)} and (A5), 
we can take the sequence $\{g_n \}_{n \in \mathbb{N}}$ such that
\begin{equation}
	g_n\in H^1 \bigl( 0,T;C(\overline{\Omega}) \bigr) {\fukao {}\cap C(\overline{Q})},
	\quad 0<C_1 \le g_n \le C_2 \quad \textrm{in~} Q
	\label{gn}
\end{equation}
for all $n \in \mathbb{N}$
and 
\begin{equation}
	g_n \to g \quad \textrm{in~} C(\overline{Q})
	\label{gn2}
\end{equation}
as $n \to \infty $, i.e.\ 
$g_n$ satisfies assumptions (A4) and (A5).

From the definition of $K(0)$, we see that 
there exists $\tilde{\sigma }_0 \in \tilde{K}(0)$ such that 
$\sigma _0=\tilde{\sigma }_0-\sigma _*(0)$. 
We define $\tilde{\sigma}_{0,n}$ by 
\begin{equation*} 
	\tilde{\sigma}_{0,n}(x):=\left(1-\frac{1}{C_1}
	\bigl| g_n(0)-g(0) \bigr|_{C(\overline{\Omega})} \right) \tilde{\sigma}_0(x)
\end{equation*}
for a.a.\ $x \in \Omega$, 
where the constant $C_1$ is the same as in (A5).  
Because
\begin{align*}
	|\tilde{\sigma}_{0,n}^D|^2 = {} &
	\left(1-\frac{1}{C_1} \bigl|g_n(0)-g(0) 
	\bigr|_{C(\overline{\Omega})}\right)^2 
	|\tilde{\sigma}_0^D|^2\\
	\le {} & 2 \left(1-\frac{1}{C_1} \bigl| g_n(0)-g(0) \bigr|_{C(\overline{\Omega})}\right)
	g(0) \\
	\le {} 
	& 2g(0)-\frac{2g(0)}{C_1} \bigl| g_n(0)-g(0) \bigr|_{C(\overline{\Omega})} \\
	\le {} 
	& 2 g(0)-2 g_n(0)+2g_n(0) - 2 \bigl| g_n(0)-g(0) \bigr|_{C(\overline{\Omega})} \\
	\le {} & 2g_n(0)
\end{align*}
a.e.\ in $\Omega $, 
$\tilde{\sigma }_{0,n} \to \tilde{\sigma }_0$ in $\mathbb{H}$ as $n \to \infty $. 
Using these sequences, 
we define a constraint set 
$K_n(t):=\tilde{K}_n(t)-\sigma_{*}(t)$ for all $t \in [0,T]$, 
where 
\begin{equation*}
	\tilde{K}_n(t):=\left\{ \tau \in \mathbb{H} : 
	\frac{1}{2} \bigl| \tau^D(x) \bigr|^2 \le g_n(t,x) \quad \mbox{for~a.a.\ } x \in \Omega \right\}.
	\label{kn}
\end{equation*}
Therefore, defining 
{\fukao $\sigma _{0,n}:=\tilde{\sigma }_{0,n}-\sigma _{*}(0)$, indeed $\sigma_* \in C([0,T];\mathbb{V})$,} we can state the following lemma:

\paragraph{Lemma 3.2.}
\textit{For the initial data $\sigma_0$, 
there exists a sequence $\{ \sigma_{0,n} \}_{n \in \mathbb{N}}$ such that 
$\sigma_{0,n} \in K_n(0)$ for all $n\in \mathbb{N}$, 
and $\sigma_{0,n} \to \sigma_0$ in $\mathbb{H}$ 
as $n \to \infty $. } \\

Under these settings, we see that {\fukao 
$\sigma _{0,n}$, $g_n$ with $K_n(t)$ satisfy assumptions (A2), (A4), and (A5) by 
replacing $\sigma _0$, $g$ with $K(t)$}, respectively.  
Therefore, applying Theorem~2.1, we see that 
there exists a unique solution $(\boldsymbol{v}_n, \sigma_n)$ 
{\fukao of modified problem for \eqref{dleqv} and \eqref{dleqs} 
in the sense of variational inequality} satisfying 
\begin{gather}
	\boldsymbol{v}'_n(t) 
	+ 
	\partial \psi \bigl( 
	\boldsymbol{v}_n(t) 
	\bigr) 
	+ 
	E_1 \sigma_n(t) 
	= \boldsymbol{f}(t)
	\quad \mbox{in~} 
	\boldsymbol{H}, 
	\ \mbox{for~a.a.\ } t \in (0,T), 
	\label{ACP1}
	\\
	\sigma'_n(t) 
	+ \partial (\varphi +I_{K_n(t)}) \bigl( \sigma_n(t) \bigr) 
	+ E_2 \boldsymbol{v}_n(t) 
	\ni 
	h(t)\quad \mbox{in~} 
	\mathbb{H}, 
	\ \mbox{for~a.a.\ } t \in (0,T), 
	\label{ACP2}
	\\
	\boldsymbol{v}_n(0)
	=\boldsymbol{v}_{0} \quad 
	\mbox{in~} \boldsymbol{H},
	\quad 
	\sigma_n(0)=\sigma_{0,n} \quad 
	\mbox{in~} \mathbb{H}.
	\notag \label{ACP3}
\end{gather}

\paragraph{Proof of Theorem~2.2.}
Firstly, we obtain a uniform estimate for $(\boldsymbol{v}_n, \sigma_n)$. 
Multiplying \eqref{ACP1} by $\boldsymbol{v}_n(s)$ at time $t=s$, and 
using the {G}auss--{G}reen relation \eqref{gg}, we have
\begin{equation}
	\frac{1}{2} \frac{d}{ds} \bigl| \boldsymbol{v}_n(s) \bigr|_{\boldsymbol{H}}^2 
	+ \nu \bigl| \boldsymbol{v}_n (s) \bigr|_{\boldsymbol{V}}^2 
	= - \bigl( \sigma_n(s), \varepsilon \bigl( \boldsymbol{v}_n(s) \bigr) \bigr)_{\mathbb{H}}
	+ \bigl( \boldsymbol{f}(s), \boldsymbol{v}_n(s) \bigr) _{\boldsymbol{H}}
	\label{aaa}
\end{equation}
for a.a.\ $s \in (0,T)$. 
In addition, from \eqref{ACP2} and \eqref{chara}, 
we have the following 
inequality:
\begin{align}
	\bigl( 
	\sigma'_n(s),\sigma_n(s)-\tau
	\bigr)_{\mathbb{H}}
	& +
	\kappa 
	\bigl( 
	\sigma_n(s),\sigma_n(s)-\tau
	\bigr)_{\mathbb{V}}
	\notag \\
	& {}-
	\bigl( 
	\varepsilon
	\bigl(
	\boldsymbol{v}_n(s)
	\bigr),
	\sigma_n(s)-\tau 
	\bigr)_{\mathbb{H}} \le 
	\bigl(
	h(s),\sigma_n(s)-\tau 
	\bigr)_{\mathbb{H}}
	\label{sub2}
\end{align}
for all $\tau \in{\risei K_n(s)}$ and for a.a.\ $s \in (0,T)$.
Now, taking $\tau :=-\sigma_*(s)$, 
\begin{align}
	& \frac{1}{2} \frac{d}{ds} \bigl| \sigma_n (s)+\sigma_{*}(s) \bigr|^2_{\mathbb{H}} 
	+ \kappa \bigl| \sigma_n (s) \bigr|_{\mathbb{V}}^2 
	\notag \\
	& \le 
	-\bigl( \sigma_{*}'(s),\sigma_n (s)+\sigma_{*}(s) \bigr)_{\mathbb{H}}
	+
	\bigl( \varepsilon \bigl( \boldsymbol{v}_n (s) \bigr), 
	\sigma_n (s)+\sigma_{*}(s) \bigr)_{\mathbb{H}} 
	+ \bigl( 
	h(s), \sigma_n (s)+\sigma_{*}(s)
	\bigr)_{\mathbb{H}} 
	\notag 
	\\
	& \quad {} - \kappa \bigl( \sigma_n(s),\sigma_{*}(s) \bigr)_{\mathbb{V}}
	\label{bbb}
\end{align}  
for a.a.\ $s \in (0,T)$. Adding these expressions 
{\fukao \eqref{aaa} and \eqref{bbb},}
and using the {Y}oung inequality, we deduce that
\begin{align*}
	& \frac{1}{2} \frac{d}{ds} \bigl| \boldsymbol{v}_n(s) \bigr|_{\boldsymbol{H}}^2 
	+\frac{1}{2} \frac{d}{ds} \bigl| \sigma_n (s)+\sigma_{*}(s) \bigr|^2_{\mathbb{H}} 
	+ \frac{\nu}{2} \bigl| \boldsymbol{v}_n (s) \bigr|_{\boldsymbol{V}}^2 
	+ \frac{\kappa}{2} \bigl| \sigma_n(s) \bigr|_{\mathbb{V}}^2
	\notag \\
	& \le \frac{1}{2} \bigl| \boldsymbol{f}(s) \bigr| _{\boldsymbol{H}}^2
	+ \frac{1}{2} \bigl| \boldsymbol{v}_n(s) \bigr| _{\boldsymbol{H}}^2 
	+ \bigl| \sigma_{*}'(s) \bigr|_{\mathbb{H}}^2
	+ \bigl| h(s) \bigr|_{\mathbb{H}}^2
	+ \frac{1}{2\nu } \bigl| \sigma_{*}(s) \bigr|_{\mathbb{H}}^2
	+ \frac{\kappa }{2} \bigl| \sigma_* (s) \bigr|_{\mathbb{V}}^2 
	\\
	& {} \quad + \frac{1}{2}\bigl| \sigma_n (s)+\sigma_{*}(s) \bigr|_{\mathbb{H}}^2 .
\end{align*}
By the {G}ronwall inequality, 
\begin{align*}
	&  
	\bigl| 
	\boldsymbol{v}_n(t)
	\bigr|^2_{\boldsymbol{H}} 
	+ \bigl| \sigma_n (t)+\sigma_{*}(t) \bigr|^2_{\mathbb{H}} 
	\notag \\
	& \le \biggl( | 
	\boldsymbol{v}_0
	|^2_{\boldsymbol{H}} 
	+ {\risei\bigl| \sigma_{0}
	+  \sigma_{*}(0) \bigr|^2_{\mathbb{H}} }
	+ | \boldsymbol{f} | _{L^2(0,T;\boldsymbol{H})}^2
	+ { \left(2+\frac{1}{\nu }\right)}| \sigma_{*}|_{H^1(0,T;\mathbb{H})}^2  
	+ { 2}| h |_{L^2(0,T;\mathbb{H})}^2 
	\notag \\
	& {} \quad + |\sigma _*|^{ 2}_{L^2(0,T;\mathbb{V})} \biggr) e^{T}
\end{align*}
for all $t \in [0,T]$.
Thus, there exists a positive constant $M_5(\nu )$ that depends on $\nu >0$, but is independent of $n \in \mathbb{N}$, such that 
\begin{equation*} 
	\bigl| 
	\boldsymbol{v}_n(t)
	\bigr|^2_{\boldsymbol{H}} 
	+\bigl| \sigma_n (t) \bigr|^2_{\mathbb{H}} 
	+ \nu \int_{0}^{t} \bigl| \boldsymbol{v}_n(s) \bigr|^2_{\boldsymbol{V}} ds
	+ \kappa \int_{0}^{t} \bigl| \sigma _n(s)\bigr|^2_{\mathbb{V}} ds
	\le M_5(\nu )
	\label{funi1}
\end{equation*}
for all $t \in [0,T]$.

Next, multiplying \eqref{ACP1} by $\boldsymbol{v}_n'(s)$ at time $t=s$, and 
using the {G}auss--{G}reen relation \eqref{gg}, 
we have
\begin{align*}
	\bigl| \boldsymbol{v}_n'(s) \bigr|_{\boldsymbol{H}}^2 
	+ \frac{\nu }{2} \frac{d}{ds} \bigl| \boldsymbol{v}_n (s) \bigr|_{\boldsymbol{V}}^2 
	= {} & \bigl( \mathbf{div}\sigma_n(s), \boldsymbol{v}_n'(s) \bigr)_{\boldsymbol{H}}
	+ \bigl( \boldsymbol{f}(s), { \boldsymbol{v}'_n(s)} \bigr) _{\boldsymbol{H}} \\
	\le {} & {\fukao 
	\bigl| \sigma _n(s) \bigr|_{\mathbb{V}}^2 
	+ \bigl| \boldsymbol{f}(s) \bigr|_{\boldsymbol{H}}^2 
	+ \frac{1}{2} \bigl| \boldsymbol{v}_n'(s) \bigr|_{\boldsymbol{H}}^2 
	}  
\end{align*}
for a.a.\ $s \in (0,T)$. Integrating this expression over $[0,t]$ with respect to time, 
we obtain 
\begin{align*}
	\int_{0}^{t}\bigl| \boldsymbol{v}_n'(s) \bigr|_{\boldsymbol{H}}^2 ds
	+ \nu  \bigl| \boldsymbol{v}_n (t) \bigr|_{\boldsymbol{V}}^2 
	\le {} &{  
	| \boldsymbol{v}_0 |_{\boldsymbol{V}}^2 
	+ \frac{2}{\kappa } M_5(\nu )
	+ 2| \boldsymbol{f} |_{L^2(0,T;\boldsymbol{H})}^2 
	}:=M_6(\kappa ,\nu )
\end{align*}
for all $t \in [0,T]$. 
By comparison, in \eqref{ACP1}, we also have
\begin{equation*}
	\int_{0}^{t}\bigl| \partial \psi \bigl( \boldsymbol{v}_n(s) \bigr) 
	\bigr|_{\boldsymbol{H}}^2 ds
	\le M_7(\kappa ,\nu )
\end{equation*}
for all $t \in [0,T]$ for some suitable positive constant $M_7(\kappa, \nu )$. 
Thus, there exists a subsequence 
$\{n_k \}_{k \in \mathbb{N}}$ with $n_k \to \infty $ as $k \to \infty $, 
$\boldsymbol{v} \in H^1(0,T;\boldsymbol{H}) \cap L^\infty (0,T;\boldsymbol{V}) \cap 
L^2(0,T;H^2(\Omega )^3)$, and $\sigma \in L^\infty (0,T;\mathbb{H}) \cap L^2(0,T;\mathbb{V})$ such that 
\begin{gather*}
	\boldsymbol{v}_{n_k} \to \boldsymbol{v} 
	\quad \textrm{weakly~star~in~} L^\infty (0,T;\boldsymbol{V}) \cap L^2\bigl(0,T;H^2(\Omega )^3 \bigr), \\
	\boldsymbol{v}_{n_k}' \to \boldsymbol{v}' 
	\quad \textrm{weakly~in~} L^2 (0,T;\boldsymbol{H}), \\
	\sigma _{n_k} \to \sigma 
	\quad \textrm{weakly~{\fukao star}~in~} L^\infty (0,T;\mathbb{H}) \cap L^2(0,T;\mathbb{V})
\end{gather*}
as $k \to \infty $. 
Moreover, from the {A}ubin compactness theorem 
(see, e.g.\ \cite{Sim87}), 
\begin{equation}
	\boldsymbol{v}_{n_k} \to \boldsymbol{v} 
	\quad \textrm{in~} C \bigl( [0,T];\boldsymbol{H} \bigr) \cap L^2 \bigl(0,T;\boldsymbol{V} \bigr)
	\label{aubin}
\end{equation}
as $k \to \infty $.

Secondly, 
to obtain strong convergence, we apply the abstract technique 
to the time-dependent constraint \cite{FK13b}. 
Recall that, for each $n \in \mathbb{N}$, there exists 
$\tilde{\sigma }_n(t) \in \tilde{K}_n(t)$
such that $\sigma_n(t)=\tilde{\sigma }_n(t)-\sigma _*(t)$
for all $t \in [0,T]$. 
For each $r\in(0,1)$, there exists $N_r \in \mathbb{N}$ such that  
\begin{equation*}
	r \tilde{\sigma}_{n_k}(t) \in \tilde{K}_{n_l}(t),
	\quad r \tilde{\sigma}_{n_l}(t)\in \tilde{K}_{n_k}(t) 
	\label{const}
\end{equation*}
for all $k,l \ge N_r$ and $t\in[0,T]$. 
Indeed, from \eqref{gn}, we see that 
$g_{n_k}(t,x)/C_1 \ge 1$ for all $(t,x) \in \overline{Q}$. 
Moreover, 
from \eqref{gn2}, there exists $N_r \in \mathbb{N}$ such that 
\begin{equation*}
	|g_{n_k}-g_{n_l}|_{C(\overline{Q})}\le C_1(1-r)
\end{equation*}
for all $k,l\ge N_r$. 
Therefore, 
\begin{align*}
	\frac{1}{2} \bigl| r \tilde{\sigma }_{n_k}^D(t) \bigr|^2 
	\le {} & \frac{1}{2} \left(1-\frac{|g_{n_k}-g_{n_l}|_{C(\overline{Q})}}{C_1}\right) 
	\bigl| \tilde{\sigma }_{n_k}^D(t) \bigr|^2 \\
	\le {} & \left(1-\frac{|g_{n_k}-g_{n_l}|_{C(\overline{Q})}}{C_1}\right)g_{n_k}(t)\\
	\le {} & g_{n_k}(t)-g_{n_l}(t)+g_{n_l}(t)-|g_{n_k}-g_{n_l}|_{C(\overline{Q})}\\
	\le {} & g_{n_l}(t)
\end{align*}
a.e.\ in $\Omega $ for all $k,l\ge N_r$. 
Therefore, $r \tilde{\sigma }_{n_k}(t) \in \tilde{K}_{n_l}(t)$ and, 
similarly, $r \tilde{\sigma }_{n_l}(t)\in \tilde{K}_{n_k}(t)$ for all $k,l\ge N_r$. 
Using this fact, we now show that 
$\{\sigma_{n_k}\}_{k\in \mathbb{N}}$ 
is a {C}auchy sequence in $C([0,T];\mathbb{H})$.
We know that $\sigma_{n_k}(s):=\tilde{\sigma }_{n_k}(s)-\sigma _*(s) \in K_{n_k}(s)$ satisfies 
\eqref{ACP2} at {\fukao time $t=s$}. Moreover, $r \tilde{\sigma }_{n_l}(s)-\sigma _*(s) \in K_{n_k}(s)$ 
for all $k,l\ge N_r$. 
Therefore, 
for all $k,l\ge N_r$, 
{\fukao taking
$\tau := r \tilde{\sigma }_{n_l}(s)-\sigma _*(s)$ as the 
test function in \eqref{sub2} of $\sigma_{n_k}(s)$}, and using \eqref{chara}, 
we obtain the following inequalities:
\begin{align}
	& \bigl( \sigma'_{n_k}(s), 
	\tilde{\sigma}_{n_k}(s)-r \tilde{\sigma}_{n_l}(s)
	\bigr)_{\mathbb{H}}
	+ \kappa \bigl( \sigma_{n_k}(s), \tilde{\sigma}_{n_k}(s)-r \tilde{\sigma}_{n_l}(s) \bigr)_{\mathbb{V}}
	\notag\\
	& \le \bigl(
	\varepsilon \bigl( \boldsymbol{v}_{n_k}(s) \bigr), 
	\tilde{\sigma}_{n_k}(s)-r\tilde{\sigma}_{n_l}(s)
	 \bigr)_{\mathbb{H}} 
	+ \bigl( h(s),\tilde{\sigma}_{n_k}(s)-r \tilde{\sigma}_{n_l}(s) \bigr)_{\mathbb{H}},
	\label{ineqk}
\end{align}
and analogously 
\begin{align}
	& \bigl( \sigma'_{n_l}(s), 
	\tilde{\sigma}_{n_l}(s)-r \tilde{\sigma}_{n_k}(s)
	\bigr)_{\mathbb{H}}
	+ \kappa \bigl( \sigma_{n_l}(s), \tilde{\sigma}_{n_l}(s)-r \tilde{\sigma}_{n_k}(s) \bigr)_{\mathbb{V}} \notag\\
	& \le \bigl(
	\varepsilon \bigl( \boldsymbol{v}_{n_l}(s) \bigr), 
	\tilde{\sigma}_{n_l}(s)-r\tilde{\sigma}_{n_k}(s)
	 \bigr)_{\mathbb{H}} 
	+ \bigl( h(s),\tilde{\sigma}_{n_l}(s)-r \tilde{\sigma}_{n_k}(s) \bigr)_{\mathbb{H}}.
	\label{ineql}
\end{align}
for a.a.\ $s \in (0,T)$. 
Taking the sum of \eqref{ineqk} and \eqref{ineql},
\begin{align}
	& \frac{1}{2} \frac{d}{ds} \bigl| \sigma_{n_k}(s)-\sigma_{n_l}(s) \bigr|^2_{\mathbb{H}} 
	+ \kappa \bigl| \sigma_{n_k}(s)-\sigma_{n_l}(s) \bigr|_{\mathbb{V}}^2
	\notag 
	\\
	& \le 
	-(1-r) \bigl( \sigma'_{n_k}(s), 
	{\fukao \sigma_{n_l}(s)+\sigma _*(s)} 
	\bigr)_{\mathbb{H}} 
	-
	(1-r) \bigl( \sigma'_{n_l}(s), 
	{\fukao \sigma_{n_k}(s)+\sigma _*(s)}
	\bigr)_{\mathbb{H}} 
	\notag 
	\\
	& \quad 
	{\fukao {} + 
	\kappa (1-r) \bigl| \sigma_{n_k}(s) \bigl|_{\mathbb{V}} 
	\bigl| \sigma_{n_l}(s) + \sigma _*(s) \bigr|_{\mathbb{V}} 
	+
	\kappa (1-r) \bigl| \sigma_{n_l}(s) \bigl|_{\mathbb{V}} 
	\bigl| \sigma_{n_\kappa }(s)+\sigma_*(s) \bigr|_{\mathbb{V}} } 
	\notag 
	\\
	& \quad {}
	+ \bigl| \boldsymbol{v}_{n_k} (s) -\boldsymbol{v}_{n_l}(s) \bigr|_{\boldsymbol{V}} 
	\bigl| \sigma _{n_k}(s)-\sigma _{n_l}(s) \bigr|_{\mathbb{H}} 
	+ (1-r) \bigl| \boldsymbol{v}_{n_k} (s) \bigr|_{\boldsymbol{V}} 
	\bigl| \sigma _{n_l}(s){\fukao {}+\sigma _*(s)} \bigr|_{\mathbb{H}} 
	\notag \\
	& \quad {} + (1-r) \bigl| \boldsymbol{v}_{n_l} (s) \bigr|_{\boldsymbol{V}} 
	\bigl| \sigma _{n_k}(s) {\fukao {}+\sigma _*(s)} \bigr|_{\mathbb{H}}
	+ (1-r) \bigl| h(t) \bigr|_{\mathbb{H}} \bigl| \sigma_{n_k}(t){\fukao {}+\sigma _*(s)} \bigr|_{\mathbb{H}}
	\notag 
	\\
	& \quad {} 
	+ (1-r) \bigl| h(t) \bigr|_{\mathbb{H}} \bigl| \sigma_{n_l}(t){\fukao {}+\sigma _*(s)} \bigr|_{\mathbb{H}},
	\label{lastesti}
\end{align}
for a.a.\ $s \in (0,T)$. 
By virtue of the {G}ronwall inequality, we deduce that 
\begin{align*}
	& \bigl| \sigma_{n_k}(t)-\sigma_{n_l}(t) \bigr|^2_{\mathbb{H}} 
	\notag 
	\\
	& \le \biggl\{ \bigl| \sigma_{0,n_k}-\sigma_{0,n_l} \bigr|^2_{\mathbb{H}} 
	+
	2(1-r) | \sigma_{0,n_l} |  
	| \sigma_{0,n_k}|_{\mathbb{H}} 
	+
	2(1-r) \bigl| \sigma_{n_k}(T) \bigr|_{\mathbb{H}} 
	\bigl| \sigma_{n_l}(T) \bigr|_{\mathbb{H}}
	\notag 
	\\ 
	& \ {\fukao {} + 2(1-r)\int_{0}^{T} 
	\bigl| \sigma_{n_k}(s) \bigl|_{\mathbb{H}} 
	\bigl| \sigma_*'(s) \bigr|_{\mathbb{H}} ds
	+ 
	2(1-r) | \sigma_{0,n_k} |  
	| \sigma_*(0)|_{\mathbb{H}} 
	+
	2(1-r) \bigl| \sigma_{n_k}(T) \bigr|_{\mathbb{H}} 
	\bigl| \sigma_*(T) \bigr|_{\mathbb{H}} }
	\notag \\
	& \ {\fukao {} + 2(1-r)\int_{0}^{T} 
	\bigl| \sigma_{n_l}(s) \bigl|_{\mathbb{H}} 
	\bigl| \sigma_*'(s) \bigr|_{\mathbb{H}} ds
	+
	2(1-r) | \sigma_{0,n_l} |  
	| \sigma_*(0)|_{\mathbb{H}} 
	+
	2(1-r) \bigl| \sigma_{n_l}(T) \bigr|_{\mathbb{H}} 
	\bigl| \sigma_*(T) \bigr|_{\mathbb{H}} }
	\notag \\
	& \ {} 
	+
	4\kappa (1-r) \int_{0}^{T}\bigl| \sigma_{n_k}(s) \bigl|_{\mathbb{V}} 
	\bigl| \sigma_{n_l}(s) \bigr|_{\mathbb{V}} ds
	+ {\fukao 2 \kappa (1-r) \int_{0}^{T}\bigl| \sigma_{n_k}(s) \bigl|_{\mathbb{V}} 
	\bigl| \sigma_*(s) \bigr|_{\mathbb{V}} ds}
	\notag \\
	& \ {} 
	+ {\fukao 2 \kappa (1-r) \int_{0}^{T}\bigl| \sigma_{n_l}(s) \bigl|_{\mathbb{V}} 
	\bigl| \sigma_*(s) \bigr|_{\mathbb{V}} ds}
	+ \int_{0}^{T} 
	\bigl| \boldsymbol{v}_{n_k} (s) -\boldsymbol{v}_{n_l}(s) \bigr|_{\boldsymbol{V}}^2 
	ds
	\notag \\
	& \ {}
	+ 2(1-r) \int_{0}^{T} \bigl| \boldsymbol{v}_{n_k} (s) \bigr|_{\boldsymbol{V}} 
	\bigl| \sigma _{n_l}(s) \bigr|_{\mathbb{H}} ds 
	{\fukao {}
	+ 2(1-r) \int_{0}^{T} \bigl| \boldsymbol{v}_{n_k} (s) \bigr|_{\boldsymbol{V}} 
	\bigl| \sigma _*(s) \bigr|_{\mathbb{H}} ds }
	\notag \\
	& \ {} 
	+ 2(1-r) \int_{0}^{T} \bigl| \boldsymbol{v}_{n_l} (s) \bigr|_{\boldsymbol{V}} 
	\bigl| \sigma _{n_k}(s) \bigr|_{\mathbb{H}} ds 
	{\fukao {}
	+ 2(1-r) \int_{0}^{T} \bigl| \boldsymbol{v}_{n_l} (s) \bigr|_{\boldsymbol{V}} 
	\bigl| \sigma _*(s) \bigr|_{\mathbb{H}} ds }
	\notag \\
	& \ {} 
	+ 2(1-r) \int_{0}^{T} 
	\bigl| h({\fukao s}) \bigr|_{\mathbb{H}} \bigl| \sigma_{n_k}({\fukao s}) \bigr|_{\mathbb{H}} ds
	{\fukao {}
	+ 2(1-r) \int_{0}^{T} 
	\bigl| h(s) \bigr|_{\mathbb{H}} \bigl| \sigma_*(s) \bigr|_{\mathbb{H}} ds }
	\notag \\
	& \ {} 
	+ 2(1-r) \int_{0}^{T} 
	\bigl| h({\fukao s}) \bigr|_{\mathbb{H}} \bigl| \sigma_{n_l}({\fukao s}) \bigr|_{\mathbb{H}} ds
	{\fukao {}
	+ 2(1-r) \int_{0}^{T} \bigl| h({\fukao s}) \bigr|_{\mathbb{H}} 
	\bigl| \sigma_*({\fukao s}) \bigr|_{\mathbb{H}} ds }
	\biggr\} e^T
\end{align*}
for all $t \in [0,T]$.
Thus, {\fukao by using \eqref{aubin} there exists a positive constant $M^*(\nu )$ such that }
\begin{equation*}
	\limsup _{k,l \to \infty } \bigl| \sigma_{n_k}(t)-\sigma_{n_l}(t) \bigr|^2_{\mathbb{H}} 
	\le (1-r) M^*(\nu ) e^T
\end{equation*}
for all $t \in [0,T]$. Letting $r \to 1$, we see that 
\begin{equation*} 
	\limsup _{k,l \to \infty } \bigl| \sigma_{n_k}(t)-\sigma_{n_l}(t) \bigr|^2_{\mathbb{H}} = 0
\end{equation*} 
for all $t \in [0,T]$, that is, $\{ \sigma_{n_k} \}_{k \in \mathbb{N}}$ is a {C}auchy sequence 
in $C([0,T];\mathbb{H})$ and 
$\sigma \in C([0,T];\mathbb{H})$. 

Next, integrating \eqref{lastesti} over $[0,T]$ with respect to time, taking 
$\limsup_{k,l \to \infty }$, and letting $r \to 1$, we finally obtain that 
\begin{equation}
	\sigma _{n_k}\to \sigma \quad \textrm{in~} C \bigl( [0,T];\mathbb{H} \bigr) \cap L^2(0,T;\mathbb{V}),
	\quad \sigma _{n_k}\to \sigma \quad \textrm{a.e.\ in~} Q
	\label{strong} 
\end{equation} 
as $k \to \infty$, and $\sigma (0)=\sigma _0$ in $\mathbb{H}$.  
Thus, setting
$\tilde{\sigma}:=\sigma+\sigma _*$, we have 
\begin{align*} 
	\frac{1}{2} \bigl| \tilde{\sigma }^D(x) \bigr|^2 
	& = \frac{1}{2} \bigl| (\sigma -\sigma _{n_k}+\sigma _{n_k}+\sigma _*)^D(x) \bigr|^2 \\
	& = \frac{1}{2} \bigl| (\sigma -\sigma _{n_k})^D(x)+ \tilde{\sigma }_{n_k}^D(x) \bigr|^2 \\
	& \le \frac{1}{2} \bigl|  (\sigma -\sigma _{n_k})^D(x) \bigl|^2 
	+ \bigl| ( \sigma -\sigma _{n_k})^D(x) \bigl| \bigl| \tilde{\sigma }_{n_k}^D(x) \bigr|
	+ \frac{1}{2} \bigl|\tilde{\sigma }_{n_k}^D(x) \bigr|^2 \\
	& \le \frac{1}{2} \bigl|  (\sigma -\sigma _{n_k})^D(x) \bigl|^2 
	+ \bigl| ( \sigma -\sigma _{n_k})^D(x) \bigl| \bigl| \tilde{\sigma }_{n_k}^D(x) \bigr|+g_{n_k}(t,x)
\end{align*} 
for a.a.\ $x \in \Omega $. Letting $k \to \infty $, we then obtain 
\begin{equation*} 
	\frac{1}{2} \bigl| \tilde{\sigma }^D(x) \bigr|^2 \le g(t,x)
\end{equation*} 
for a.a.\ $x \in \Omega $, 
that is, $\sigma (t) \in K(t)$ for a.a.\ $t \in (0,T)$.

Finally, we show that the pair $(\boldsymbol{v}, \sigma)$ 
satisfies the weak variational inequality 
\eqref{wvi2}. 
For each test function $\eta \in {\mathcal K}_0$, there exists 
a function $\tilde{\eta }:(0,T) \to \mathbb{H}$ such that 
$\tilde{\eta }(t) \in \tilde{K}(t)$ and 
$\eta(t) =\tilde{\eta }(t)-\sigma _*(t)$ for a.a.\ $t \in (0,T)$. 
Moreover, from assumption (A3), we see that 
$\tilde{\eta }=\eta +\sigma _* \in H^1(0,T;\mathbb{H})\cap L^2(0,T;\mathbb{V})$.  
Now, 
for each $r\in(0,1)$, there exists $N_r^* \in \mathbb{N}$ such that  
\begin{equation*}
	r \tilde{\eta }(t) \in \tilde{K}_{n_k}(t) 
\end{equation*}
for all $k \ge N_r^*$ and $t\in[0,T]$. 
Indeed, from \eqref{gn}, we see that 
$g_{n_k}(t,x)/C_1 \ge 1$ for all $(t,x) \in \overline{Q}$. 
Moreover, 
from \eqref{gn2} there exists $N_r^* \in \mathbb{N}$ such that 
\begin{equation*}
	|g_{n_k}-g|_{C(\overline{Q})}\le C_1(1-r)
\end{equation*}
for all $k \ge N_r^*$. 
Therefore, 
\begin{align*}
	\frac{1}{2} \bigl| r \tilde{\eta}^D(t) \bigr|^2 
	\le {} & \frac{1}{2} \left(1-\frac{|g_{n_k}-g|_{C(\overline{Q})}}{C_1}\right) 
	\bigl| \tilde{\eta }^D(t) \bigr|^2 \\
	\le {} & \left(1-\frac{|g_{n_k}-g|_{C(\overline{Q})}}{C_1}\right)g(t)\\
	\le {} & g(t)-g_{n_k}(t)+g_{n_k}(t)-|g_{n_k}-g|_{C(\overline{Q})}\\
	\le {} & g_{n_k}(t)
\end{align*}
a.e.\ in $\Omega $ for all $k \ge N_r$, namely 
$r \tilde{\eta}(t) \in \tilde{K}_{n_k}(t)$ for all $k \ge N_r$. 
Now, for all $k \ge N_r$, 
{\fukao taking
$\tau := r \tilde{\eta}(s)-\sigma _*(s)\in K_{n_k}(s)$ as the 
test function in \eqref{sub2} of $\sigma_{n_k}(s)$}
and integrating the resultant over $(0,t)$ with respect to time, we obtain 
\begin{align}
	& \int_{0}^{t} \bigl( 
	\sigma_{n_k}'(s),  \tilde{\sigma}_{n_k}(s)-\tilde{\eta}(s)
	\bigr)_{\mathbb{H}} ds 
	+ (1-r)\int_{0}^{t} \bigl( 
	\sigma_{n_k}'(s), \tilde{\eta}(s)
	\bigr)_{\mathbb{H}} ds 
	\notag \\
	& {}+\kappa \int_{0}^{t} \bigl( 
	\sigma_{n_k}(s), \tilde{\sigma}_{n_k}(s)-\tilde{\eta}(s)
	\bigr)_{\mathbb{V}} ds 
	+(1-r)\kappa \int_{0}^{t} \bigl( 
	\sigma_{n_k}(s), \tilde{\eta}(s)
	\bigr)_{\mathbb{V}} ds 
	\notag \\
	& {} - 
	\int_{0}^{t}
	\bigl( 
	\varepsilon \bigl( 
	\boldsymbol{v}_{n_k}(s)
	\bigr), 
	\tilde{\sigma}_{n_k}(s)- \tilde{\eta}(s)
	\bigr)_{\mathbb{H}} ds 
	- (1-r)
	\int_{0}^{t}
	\bigl( 
	\varepsilon \bigl( 
	\boldsymbol{v}_{n_k}(s)
	\bigr), 
	\tilde{\eta}(s)
	\bigr)_{\mathbb{H}} ds 
	\notag \\
	& \quad {}
	\le \int_{0}^{t}
	\bigl( 
	h(s), 
	\tilde{\sigma}_{n_k}(s) - \tilde{\eta}(s)
	\bigr)_{\mathbb{H}} ds 
	+ (1-r)
	\int_{0}^{t}
	\bigl( 
	h(s), 
	\tilde{\sigma}_{n_k}(s) - \tilde{\eta}(s)
	\bigr)_{\mathbb{H}} ds
	\label{wvi2n}
\end{align}
for all $t \in [0,T]$.
Here, as we know that 
\begin{align*}
	& \int_{0}^{t} \bigl( 
	\sigma_{n_k}'(s),  \tilde{\sigma}_{n_k}(s)-\tilde{\eta}(s)
	\bigr)_{\mathbb{H}} ds 
	\notag \\
	& = \int_{0}^{t} \bigl( 
	\sigma_{n_k}'(s)-\eta' (s), \sigma_{n_k}(s)-\eta(s)
	\bigr)_{\mathbb{H}} ds 
	+\int_{0}^{t} \bigl( 
	\eta' (s), \sigma_{n_k}(s)-\eta(s)
	\bigr)_{\mathbb{H}} ds 
	\notag \\
	& = 
	\frac{1}{2} \bigl| \sigma_{n_k} (t) - \eta (t) \bigr|_{\mathbb{H}}^2 
	-\frac{1}{2}  \bigl| \sigma_{\risei{0,n_k}} - \eta (0) \bigr|_{\mathbb{H}}^2 
	+ \int_{0}^{t} \bigl( 
	\eta' (s), \sigma_{n_k}(s)-\eta(s)
	\bigr)_{\mathbb{H}} ds,
\end{align*}
and
\begin{align*}
	&
	(1-r)\int_{0}^{t} \bigl( 
	\sigma_{n_k}'(s), \tilde{\eta}(s)
	\bigr)_{\mathbb{H}} ds 
	\notag \\
	&=  
	(1-r)\bigl( 
	\sigma_{n_k}(t), \tilde{\eta}(t)
	\bigr)_{\mathbb{H}}
	-(1-r)\bigl( 
	\sigma_{\risei{0,n_k}}, \tilde{\eta}(0)
	\bigr)_{\mathbb{H}}
	+(1-r)\int_{0}^{t} \bigl( 
	\sigma_{n_k}(s), \tilde{\eta}'(s)
	\bigr)_{\mathbb{H}} ds
\end{align*}
for all $t \in [0,T]$, 
we let $k \to \infty $ in \eqref{wvi2n} with the convergence
shown in {\fukao \eqref{aubin} and \eqref{strong}, use \eqref{ddd}} and deduce
\begin{align*}
	& \int_{0}^{t} \bigl( 
	\eta' (s), \sigma(s)-\eta(s)
	\bigr)_{\mathbb{H}} ds
	+\kappa \int_{0}^{t} \bigl( 
	\sigma(s), \sigma(s)- \eta(s)
	\bigr)_{\mathbb{V}} ds 
	\notag \\
	& - 
	\int_{0}^{t}
	\bigl( 
	\varepsilon \bigl( 
	\boldsymbol{v}(s)
	\bigr), 
	\sigma(s)- \eta(s)
	\bigr)_{\mathbb{H}} ds 
	+ \frac{1}{2} \bigl| \sigma (t) - \eta (t) \bigr|_{\mathbb{H}}^2 
	\notag \\
	{}
	& \quad \le \int_{0}^{t}
	\bigl( 
	h(s), 
	\sigma(s)- \eta(s)
	\bigr)_{\mathbb{H}} ds 
	+ \frac{1}{2}  \bigl| \sigma_0 - \eta (0) \bigr|_{\mathbb{H}}^2 
	\notag \\
	& \quad \quad {}
	-(1-r)\bigl( 
	\sigma(t), \tilde{\eta}(t)
	\bigr)_{\mathbb{H}}
	+(1-r)\bigl( 
	\sigma_{0}, \tilde{\eta}(0)
	\bigr)_{\mathbb{H}}
	-{\fukao (1-r)} \int_{0}^{t} \bigl( 
	\sigma(s), \tilde{\eta}'(s)
	\bigr)_{\mathbb{H}} ds
	\notag \\
	& \quad \quad {}
	-(1-r)\kappa \int_{0}^{t} \bigl( 
	\sigma(s), \tilde{\eta}(s)
	\bigr)_{\mathbb{V}} ds 
	+(1-r)
	\int_{0}^{t}
	\bigl( 
	\varepsilon \bigl( 
	\boldsymbol{v}(s)
	\bigr), 
	\tilde{\eta}(s)
	\bigr)_{\mathbb{H}} ds 
	\notag \\
	& \quad \quad {} 
	+ (1-r)
	\int_{0}^{t}
	\bigl( 
	h(s), 
	\sigma(s) - \eta(s)
	\bigr)_{\mathbb{H}} ds
\end{align*}
for all $t \in [0,T]$. Thus, letting $r \to 1$, we see that $(\boldsymbol{v}, \sigma )$ satisfies 
the weak variational inequality \eqref{wvi2}. 
It is clear that $(\boldsymbol{v}, \sigma )$ satisfies \eqref{wvi1}. Therefore, 
$(\boldsymbol{v}, \sigma )$ is a solution 
{\fukao of modified problem for \eqref{dleqv} and \eqref{dleqs}}
in the sense of weak variational inequality.

The proof of the 
uniqueness in the sense of weak variational inequality is the same as 
in \cite[Theorem~2.1]{FK13b}. Therefore, we omit it from this paper. 
\hfill $\Box$

\section{Viscous perfect plasticity model}

In the final section, we discuss the existence problem for the {\fukao viscous} perfect plasticity model.

\subsection{Auxiliary problems}

Conisider two auxiliary problems for the {\fukao viscous} perfect plasticity model. 
{\fukao We} obtain the following lemma.

\paragraph{Lemma 4.1.}
\textit{The set of proper, lower semicontinuous, 
and convex functions 
$\{ I_{K(t)} \}$ satisfies condition {\upshape (H)} {\fukao for all $s, t \in (0,T)$ with $|t-s|<T_*$}. }\\

{\fukao 
\paragraph{Proof.}
It is clear that $I_{K(t)}$ 
is a lower semicontinuous functional on $\mathbb{H}$.
We show that 
$\{ I_{K(t)} \}$ satisfies condition {\upshape (H)}. 
The proof is essentially same as one of Lemma~3.1. 
For each $r>0$, $s,t \in [0,T]$ with $s \le t$ and $\tau \in K(s)$ with 
$|\tau |_{\mathbb{H}} \le r$, from the definition of 
$K(s)$, we see that there exists $\tilde{\tau } \in \tilde{K}(s)$ such that 
$\tau =\tilde{\tau}-\sigma _*(s)$. 
Now, we take 
\begin{equation*}
	\tilde{\tau}_{*}:=\left(1-\frac{1}{C_1} \bigl| g(t)-g(s) 
	\bigr|_{C(\overline{\Omega})} \right) \tilde{\tau}.
\end{equation*}
Then, $\tilde{\tau }_{*} \in \tilde{K}(t)$ from the same proof of Lemma~3.1, that is, 
$\tau_{*}:=\tilde{\tau}_{*}-\sigma_{*}(t)$ is an element of $K(t)$. 
We obtain
\begin{align*}
	|\tau_{*}-\tau|_{\mathbb{H}} 
	& \le 
	|\tilde{\tau}_{*}-\tilde{\tau}|_{\mathbb{H}} 
	+ 
	\bigl| \sigma_{*}(t)-\sigma_{*}(s) \bigr|_{\mathbb{H}}
	\\
	& \le \int_s^t \left( \frac{1}{C_1}
	\bigl| g'(l) \bigr|_{C(\overline{\Omega})} \bigl| \tilde{\tau }-\sigma _*(s) \bigr|_{\mathbb{H}}
	+\frac{1}{C_1}
	\bigl| g'(l) \bigr|_{C(\overline{\Omega})} |\sigma _*(s)|_{\mathbb{H}}
	+
	\bigl| 
	\sigma_{*}'(l) 
	\bigr|_{\mathbb{H}} 
	\right) dl 
	\\
	& \le \int_{s}^{t} \alpha _r (l) dl 
\end{align*}
for all $s, t \in (0,T)$ with $|t-s|<T_*$, where $T_*$ is same as in Lemma~3.1. 
Therefore, we can 
take function $\alpha_r(\cdot ):=
(1/C_1)|g'(\cdot )|_{C(\overline{\Omega})}(r+|\sigma _*|_{C([0,T];\mathbb{H})})
+|\sigma_{*}'(\cdot )|_{\mathbb{H}}$ for each $r>0$. 

Next, $I_{K(s)}(\tau )=I_{K(t)}(\tau _*)=0$. Therefore, 
taking $\beta_r(\cdot ):=0$ we see that 
\eqref{H2} holds for all $s, t \in (0,T)$ with $|t-s|<T_*$. } \hfill $\Box$ \\

Using Lemma~4.1, we can apply Proposition~3.1 again to 
obtain the solution $\sigma \in H^1(0,T;\mathbb{H})$ of the following 
form.

\paragraph{Proposition 4.1.} 
\textit{For any given 
$\tilde{\boldsymbol{v}} 
\in L^2(0,T;\boldsymbol{V})$, 
there exists a unique solution 
$\sigma \in H^1(0,T;\mathbb{H})$, with 
$\sigma (t) \in K(t)$ for all $t \in [0,T]$, 
to the following problem.
\begin{gather}
	\sigma'(t) 
	+ \partial I_{K(t)}
	\bigl( \sigma(t) \bigr) + 
	E_2 \tilde{\boldsymbol{v}}(t) 
	\ni h(t) 
	\quad \mbox{in~} \mathbb{H}, 
	\ \mbox{for~a.a.\ } t \in (0,T),
	\label{prosig12}\\
	\sigma(0)=
	\sigma_0 \quad \mbox{in~}
	\mathbb{H}.
	\label{prosig22}
\end{gather}
Moreover, there exist positive constants $M_6$ and $M_7$, 
independent of $\nu$, such that
\begin{gather} 
	\bigl| 
	\sigma(t)
	\bigr|^2_{\mathbb{H}} 
	\le M_6 \left( 1+
	\int_{0}^{T}\bigl| \tilde{\boldsymbol{v}}(s) \bigr|_{\boldsymbol{V}}^2 ds \right),
	\label{prosig32}
	\\
	\int_{0}^{t} \bigl| \sigma'(s) \bigr|_{\mathbb{H}}^2 ds
	\le  
	M_7  \left( 1+
	\int_{0}^{T}\bigl| \tilde{\boldsymbol{v}}(s) \bigr|_{\boldsymbol{V}}^2 ds \right)
	\label{prosig42}
\end{gather}
for all $t \in [0,T]$. }

\paragraph{Proof.}
The proof {\fukao is completely same as Proposition~3.2}. Here, we {\fukao only treat} the 
uniform estimate. 
By the definition of the subdifferential, \eqref{prosig12} is equivalent to the following 
inequality:
\begin{equation}
	\bigl( 
	\sigma'(s),\sigma(s)-\tau
	\bigr)_{\mathbb{H}}
	-
	\bigl( 
	\varepsilon
	\bigl(
	\tilde{\boldsymbol{v}}(s)
	\bigr),
	\sigma(s)-\tau 
	\bigr)_{\mathbb{H}} \le 
	\bigl(
	h(s),\sigma(s)-\tau 
	\bigr)_{\mathbb{H}}
	\label{sub0}
\end{equation}
for all $\tau \in K(s)$ and for a.a.\ $s \in (0,T)$. From 
assumption (A3), we can 
substitute $\tau:=-\sigma_{*}(s) \in K(s)$ in \eqref{sub0} to obtain
\begin{align*}
	& \frac{1}{2} \frac{d}{ds} \bigl| \sigma(s)+\sigma_{*}(s) \bigr|^2_{\mathbb{H}} \notag \\
	& \le 
	-\bigl( \sigma_{*}'(s),\sigma(s)+\sigma_{*}(s) \bigr)_{\mathbb{H}}
	+
	\bigl( \varepsilon \bigl( \tilde{\boldsymbol{v}} (s) \bigr), 
	\sigma(s)+\sigma_{*}(s) \bigr)_{\mathbb{H}} 
	+ \bigl( 
	h(s), \sigma(s)+\sigma_{*}(s)
	\bigr)_{\mathbb{H}} \notag \\
	& \le
	\frac{3}{2} \bigl| \sigma _{*}'(s) \bigr|_{\mathbb{H}}^2
	+
	\frac{3}{2} \bigl| \tilde{\boldsymbol{v}}(s) \bigr|_{\boldsymbol{V}}^2  
	+ \frac{3}{2} \bigl| h(s) \bigr|_{\mathbb{H}}^2 + \frac{1}{2} \bigl| 
	\sigma(s)+\sigma_{*}(s) \bigr|_{\mathbb{H}}^2, \label{uni12}
\end{align*}  
that is,
\begin{equation*}
	\frac{d}{ds} \bigl| 
	\sigma(s)+\sigma_{*}(s)
	\bigr|^2_{\mathbb{H}} \notag\\
	\le 
	3 \bigl| \sigma _{*}'(s) \bigr|_{\mathbb{H}}^2
	+ 3 \bigl| \tilde{\boldsymbol{v}}(s) \bigr|_{\boldsymbol{V}}^2  
	+ 3 \bigl| h(s) \bigr|_{\mathbb{H}}^2 + \bigl| 
	\sigma(s)+\sigma_{*}(s)
	\bigr|^2_{\mathbb{H}} 
\end{equation*}  
for a.a.\ $s \in (0,T)$. 
Now, using the {G}ronwall inequality, we see that 
there exists a constant $\tilde{M}_6$ that depends on 
$|\sigma_0 |_{\mathbb{H}}$, $|\sigma_*(0) |_{\mathbb{H}}$,  
$|\sigma_*'|_{L^2(0,T;\mathbb{H})}$, 
$|h|_{L^2(0,T;\mathbb{H})}$, and $T$ such that
\begin{equation*} 
	\bigl| 
	\sigma(t)+\sigma_{*}(t)
	\bigr|^2_{\mathbb{H}} 
	\le \tilde{M}_6 \left( 1+
	\int_{0}^{T}\bigl| \tilde{\boldsymbol{v}}(s) \bigr|_{\boldsymbol{V}}^2 ds \right).
\end{equation*} 
for all $t \in [0,T]$. 
This implies \eqref{prosig32} holds for some suitable constant $M_6>0$. 

Next, we recall the the approximate problem of \eqref{prosig12} using the 
{M}oreau--{Y}osida regularization $I_{K(s)}^\lambda $ for $\lambda \in (0,1]$. 
If we let $P_{K(s)}$ be the projection on $K(s)$, this is the same as the resolvent 
$(I+\lambda \partial I_{K(s)})^{-1}$. 
From Lemma~4.1, {\fukao for each $r>0$, $s,t \in [0,T]$ and 
$\tau \in \mathbb{H}$ with $|\tau|_{\mathbb{H}}\le r$}, there exists $\tilde{\tau } \in K(t)$ such that 
\begin{equation*}
	|\tilde{\tau }-P_{K(s)} \tau |_{\mathbb{H}} 
	\le \int_{s}^{t} {\fukao \alpha_r} (l) dl. 
\end{equation*}
Hence, 
\begin{align*}
	I_{K(t)}^\lambda (\tau )
	-I_{K(s)}^\lambda (\tau ) 
	\le {} & \frac{1}{2\lambda } |\tau -\tilde{\tau }|_{\mathbb{H}}^2 
	- \frac{1}{2\lambda } |\tau -P_{K(s)} \tau |_{\mathbb{H}}^2 \\
	= {} & \frac{1}{2\lambda } \bigl( |\tau -\tilde{\tau }|_{\mathbb{H}}
	+|\tau -P_{K(s)} \tau |_{\mathbb{H}} \bigr) 
	\bigl( |\tau -\tilde{\tau }|_{\mathbb{H}}
	-|\tau -P_{K(s)} \tau |_{\mathbb{H}} \bigr) 
	\\
	\le  {} & \frac{1}{2\lambda } 
	\bigl( |\tau -P_{K(s)} \tau+P_{K(s)} \tau-\tilde{\tau }|_{\mathbb{H}}
	+|\tau -P_{K(s)} \tau |_{\mathbb{H}} \bigr) 
	|\tilde{\tau }-P_{K(s)} \tau |_{\mathbb{H}}  
	\\
	\le {} & |\tilde{\tau }-P_{K(s)} \tau|_{\mathbb{H}} 
	\left| \frac{\tau -P_{K(s)} \tau }{\lambda }\right|_{\mathbb{H}}
	+ \frac{1}{2\lambda }|\tilde{\tau } -P_{K(s)} \tau |_{\mathbb{H}}^2 \\
	\le {} & \bigl| \partial I_{K(s)}^{{\fukao \lambda }} (\tau )\bigr|_{\mathbb{H}} 
	\int_{s}^{t} \alpha _{\fukao r}(l) dl 
	 + \frac{t-s}{2\lambda } \int_{s}^{t} \bigl| \alpha_{\fukao r} (l) \bigr|^2 dl. 
\end{align*}
This means that 
\begin{equation*}
	\frac{d}{ds} I_{K(s)}^\lambda (\tau ) \le 
	\bigl| \partial I_{K(s)}^{{\fukao \lambda }} (\tau )
	\bigr|_{\mathbb{H}} \bigl| \alpha_{\fukao r} (s) \bigr| 
\end{equation*}
for a.a.\ $s \in (0,T)$, namely 
\begin{equation*}
	\frac{d}{ds} I_{K(s)}^\lambda ({\fukao \eta}(s)) 
	- \bigl( {\fukao \eta}'(s), \partial I_{K(s)}^\lambda  \bigl( 
	{\fukao \eta}(s) \bigr) \bigr)_{\mathbb{H}} 
	\le 
	\bigl| \partial I_{K(s)}^\lambda  ({\fukao \eta})\bigr|_{\mathbb{H}} 
	\bigl| \alpha _{\fukao r}(s) \bigr| 
\end{equation*}
for a.a.\ $s \in (0,T)$ and for all 
${\fukao \eta} \in W^{1,1}(0,T;\mathbb{H})$. 
As a result of this rigorous estimate, we can  
multiply \eqref{prosig12} by $\partial I_{K(s)}^{{\fukao \lambda }}(\sigma_{\fukao \lambda } (s))$ at time $t=s$ 
to obtain
\begin{align*}
	& \bigl| \partial I_{K(s)}^{{\fukao \lambda }} \bigl( \sigma_{{\fukao \lambda }}(s) \bigr) \bigr|_{\mathbb{H}}^2 
	+ \frac{d}{ds} I_{K(s)}^{{\fukao \lambda }} \bigl( \sigma_{{\fukao \lambda }}(s) \bigr) \\
	& \le  
	\bigl| \partial I_{K(s)}^{{\fukao \lambda }} \bigl(\sigma_{{\fukao \lambda }}(s) 
	\bigr) \bigr|_{\mathbb{H}} \bigl| \alpha _{{\fukao r}}(s) \bigr| 
	- \bigl( \varepsilon \bigl( \tilde{\boldsymbol{v}}(s) \bigr), \partial I_{K(s)}^{{\fukao \lambda }} 
	\bigl( \sigma_{{\fukao \lambda }}(s) \bigr)
	\bigr)_{\mathbb{H}}
	+
	\bigl( h(s),  \partial I_{K(s)}^{{\fukao \lambda }} \bigl( \sigma_{{\fukao \lambda }}(s) \bigr)
	\bigr)_{\mathbb{H}} \notag \\
	& \le  
	\frac{1}{2} \bigl| \partial I_{K(s)}^{{\fukao \lambda }} 
	\bigl( \sigma_{{\fukao \lambda }}(s) \bigr) \bigr|_{\mathbb{H}}^2 +
	\frac{3}{2} \bigl| \alpha_{{\fukao r}} (s) \bigr|^2
	+ \frac{3}{2}
	\bigl| \tilde{\boldsymbol{v}}(s) \bigr|_{\boldsymbol{V}}^2 
	+ \frac{3}{2}
	\bigl|h(s) 
	\bigr|_{\mathbb{H}}^2
\end{align*}
for a.a.\ $s \in (0,T)$.
Integrating the above resultant over $[0,t]$ with respect to time, we have
\begin{align*}
	\int_{0}^{t} \bigl| \partial I_{K(s)}^{{\fukao \lambda }} 
	\bigl( \sigma_{{\fukao \lambda }}(s)\bigr)  \bigr|_{\mathbb{H}}^2 ds
	+ 2 I_{K(t)}^{{\fukao \lambda }} \bigl( \sigma_{{\fukao \lambda }}(t) \bigr) 
	\le {} & 
	3 \int_{0}^{T}\bigl| \alpha _{{\fukao r}}(s) \bigr|^2 ds
	+
	3\int_{0}^{T}
	\bigl| \tilde{\boldsymbol{v}}(s) \bigr|_{\boldsymbol{V}}^2 ds 
	+ 3 \int_{0}^{T}
	\bigl|h(s) 
	\bigr|_{\mathbb{H}}^2 ds \\
	\le {} & 
	\tilde{M}_7  \left( 1+
	\int_{0}^{T}\bigl| \tilde{\boldsymbol{v}}(s) \bigr|_{\boldsymbol{V}}^2 ds \right)
\end{align*}
for all $t \in [0,T]$, 
where $\tilde{M}_7:=\max\{3, |\alpha _{\fukao r}|_{L^2(0,T)}^2+|h|_{L^2(0,T;\mathbb{H})}^2 \}$. 
Finally, by comparison, we obtain \eqref{prosig42} with some suitable constant $M_7>0$. \hfill $\Box$ \\

We now define a solution operator $S_3:L^2(0,T;\boldsymbol{V}) \to L^2(0,T;\mathbb{H})$ 
that assigns a unique solution $S_3 \tilde{\boldsymbol{v}}:=\sigma$ 
to the above problem \eqref{prosig12}--\eqref{prosig22} on $[0,T]$.

Define $E_3:\mathbb{H} \to \boldsymbol{V}^*$ by 
\begin{equation*} 
	\langle E_3 \tau , \boldsymbol{z} \rangle _{\boldsymbol{V}^*,\boldsymbol{V}}
	:=\bigl( \varepsilon (\boldsymbol{z}), \tau \bigr)_{\mathbb{H}} 
	\quad \textrm{for~all~} {\risei\tau\in\mathbb{H},}~\boldsymbol{z} \in \boldsymbol{V}. 
\end{equation*} 
If $E_3 \tau \in \boldsymbol{H}$ then, from \eqref{gg}, we have $E_3 \tau =E_1 \tau$. 
We consider the following auxiliary problem.

\paragraph{Proposition 4.2.}
\textit{Let $\nu \in (0,1]$. 
For any given $\tilde{\sigma} \in { L^2(0,T;\mathbb{H})}$, there exists a unique solution 
$\boldsymbol{v} \in H^1(0,T;\boldsymbol{V}^*) 
\cap L^\infty (0,T;\boldsymbol{H}) \cap L^2(0,T;\boldsymbol{V})$ 
to the following problem: 
\begin{gather}
	{\fukao 
	\bigl\langle 
	\boldsymbol{v}'(t), \boldsymbol{z} 
	\bigr\rangle _{\boldsymbol{V}^*,\boldsymbol{V}} 
	+ \nu \bigl( \! \bigl( \boldsymbol{v}(t), \boldsymbol{z} \bigr) \! \bigr)
	+ \bigl\langle E_3 \tilde{\sigma}(t), \boldsymbol{z} 
	\bigr\rangle_{\boldsymbol{V}^*,\boldsymbol{V}}  
	= \bigl( \boldsymbol{f}(t), \boldsymbol{z} \bigr)_{\boldsymbol{H}} 
	} 
	\notag \\ 
	{\fukao \quad \mbox{for~all~} \boldsymbol{z} \in \boldsymbol{V}}, 
	\ \mbox{for~a.a.~} t \in (0,T), 
	\label{prov12}
	\\
	\boldsymbol{v}(0)=\boldsymbol{v}_0 \quad \mbox{in~} \boldsymbol{H}.
	\label{prov22} 
\end{gather}
Moreover, there exists a positive constant $M_8$ that is independent of $\nu $ such that
\begin{equation} 
	\bigl| 
	\boldsymbol{v}(t)
	\bigr|^2_{\boldsymbol{H}} 
	+\nu \int_{0}^{t} \bigl| \boldsymbol{v}(s) \bigr|^2_{\boldsymbol{V}} ds
	\le M_8 \left( 1+ \frac{1}{\nu }
	\int_{0}^{T} \bigl| \tilde{\sigma }(s) \bigr|_{\mathbb{H}}^2 ds \right)
	\label{prov32}
\end{equation}
for all $t \in [0,T]$. } \\

\paragraph{Proof.} 
The proof of existence {\risei and uniqueness} is quite standard. 
{\fukao Therefore, we only show the sketch of the existence proof. 
For $\boldsymbol{f}-E_3 \tilde{\sigma} \in L^2(0,T;\boldsymbol{V}^*)$, 
there exists $\{ \tilde{\boldsymbol{f}}_n \}_{n \in \mathbb{N}} \subset L^2(0,T;\boldsymbol{H})$ such that 
$\tilde{\boldsymbol{f}}_n \to \boldsymbol{f}-E_3 \tilde{\sigma}$ in 
$L^2(0,T;\boldsymbol{V}^*)$ as $n \to \infty $. Therefore, we consider the following 
approximate problem: 
\begin{gather}
	\boldsymbol{v}_n'(t) + \partial \psi \bigl( \boldsymbol{v}_n(t) \bigr) 
	= \tilde{\boldsymbol{f}}_n(t) \quad \mbox{in~} \boldsymbol{H}, \ \mbox{for~a.a.~} t \in (0,T), 
	\label{eee}
	\\
	\boldsymbol{v}_n(0)=\boldsymbol{v}_0 \quad \mbox{in~} \boldsymbol{H}.
	\label{fff}
\end{gather}
Multiplying \eqref{eee} by $\boldsymbol{v}_n(s)$ at time $t=s$, we obtain
\begin{align*}
	\frac{1}{2} \frac{d}{ds} \bigl| \boldsymbol{v}_n(s) \bigr|_{\boldsymbol{H}}^2 
	+ \nu \bigl| \boldsymbol{v}_n(s) \bigr|_{\boldsymbol{V}}^2 
	= {} &  \bigl\langle  \tilde{\boldsymbol{f}}_n(s), \boldsymbol{v}_n(s) \bigr\rangle_{\boldsymbol{V}^*,\boldsymbol{V}} \\
	\le {} & 
	\frac{1}{2\nu } \bigl| \tilde{\boldsymbol{f}}_n(s) \bigr|_{\boldsymbol{V}^*}^2 
	+ \frac{\nu }{2} \bigl| \boldsymbol{v}_n (s) \bigr|_{\boldsymbol{V}}^2 
\end{align*}
for a.a.\ $s \in (0,T)$. Therefore, using the {G}ronwall inequality, we obtain
\begin{gather*}
	\bigl| \boldsymbol{v}_n(t) \bigr|_{\boldsymbol{H}}^2 
	\le \left( 
	\bigl| \boldsymbol{v}_0 \bigr|_{\boldsymbol{H}}^2 
	+ \frac{1}{\nu } \int_{0}^{T} \bigl| \tilde{\boldsymbol{f}}_n(s) \bigr|_{\boldsymbol{V}^*}^2 ds 
	\right) 
	e^{T}, 
	\\
	\int_{0}^{t} 
	\bigl| \boldsymbol{v}_n(s) \bigr|_{\boldsymbol{V}}^2 ds 
	\le 
	\frac{1}{\nu }
	\left( 
	\bigl| \boldsymbol{v}_0 \bigr|_{\boldsymbol{H}}^2 
	+ \frac{1}{\nu } \int_{0}^{T} \bigl| \tilde{\boldsymbol{f}}_n(s) \bigr|_{\boldsymbol{V}^*}^2 ds 
	\right) 
\end{gather*}
for all $t \in [0,T]$. Moreover, from \eqref{eee} we see that 
\begin{equation*}
	\bigl\langle 
	\boldsymbol{v}'_n(t), \boldsymbol{z} 
	\bigr\rangle _{\boldsymbol{V}^*,\boldsymbol{V}} 
	+ \nu \bigl( \! \bigl( \boldsymbol{v}_n(t), \boldsymbol{z} \bigr) \! \bigr)
	= \bigl\langle \tilde{\boldsymbol{f}}_n(t), \boldsymbol{z} \bigr\rangle_{\boldsymbol{V}^*,\boldsymbol{V}}  
	\quad \mbox{for~all~} \boldsymbol{z} \in \boldsymbol{V}, 
\end{equation*}
for a.a.\ $t \in (0,T)$, namely $\{ \boldsymbol{v}'_n \}_{n \in \mathbb{N}}$ is bounded in $L^2(0,T;\boldsymbol{V}^*)$. 
Thus, there exists 
$\boldsymbol{v} \in H^1(0,T;\boldsymbol{V}^*) 
\cap L^\infty (0,T;\boldsymbol{H}) \cap L^2(0,T;\boldsymbol{V})$ such that $\boldsymbol{v}$ satisfies 
\eqref{eee} and \eqref{fff}. }

{\fukao Next, we} show the uniform estimates. 
{\fukao Taking $\boldsymbol{z}:=\boldsymbol{v}(s)$ as the test function in \eqref{prov12}} at time $t=s$, we obtain
\begin{align*}
	\frac{1}{2} \frac{d}{ds} \bigl| \boldsymbol{v}(s) \bigr|_{\boldsymbol{H}}^2 
	+ \nu \bigl| \boldsymbol{v} (s) \bigr|_{\boldsymbol{V}}^2 
	= {} & - \bigl( \tilde{\sigma }(s), \varepsilon \bigl( \boldsymbol{v}(s) \bigr) \bigr)_{\mathbb{H}}
	+ \bigl( \boldsymbol{f}(s), \boldsymbol{v}(s) \bigr) _{\boldsymbol{H}} \\
	\le {} & \frac{1}{2\nu } \bigl| \tilde{\sigma }(s) \bigr|_{\mathbb{H}}^2 
	+ \frac{\nu }{2}\bigl| \boldsymbol{v} (s) \bigr|_{\boldsymbol{V}}^2 
	+ \frac{1}{2} \bigl| \boldsymbol{f}(s) \bigr|_{\boldsymbol{H}}^2 
	+ \frac{1}{2} \bigl| \boldsymbol{v} (s) \bigr|_{\boldsymbol{H}}^2 
\end{align*}
for a.a.\ $s \in (0,T)$. Therefore, using the {G}ronwall inequality, we have
\begin{align*}
	\bigl| \boldsymbol{v}(t) \bigr|_{\boldsymbol{H}}^2 
	\le {} & \left( 
	\bigl| \boldsymbol{v}_0 \bigr|_{\boldsymbol{H}}^2 
	+ \frac{1}{\nu } \int_{0}^{T} 
	\bigl| \tilde{\sigma }(s) \bigr|_{\mathbb{H}}^2 ds 
	+ \int_{0}^{T} \bigl| \boldsymbol{f}(s) \bigr|_{\boldsymbol{H}}^2 ds 
	\right) 
	e^{T} \\
	\le {} & \tilde{M}_8 \left( 1+ \frac{1}{\nu }
	\int_{0}^{T} \bigl| \tilde{\sigma }(s) \bigr|_{\mathbb{H}}^2 ds \right)
\end{align*}
for all $t \in [0,T]$, where $\tilde{M}_8:=e^T\max\{1, 
| \boldsymbol{v}_0 |_{\boldsymbol{H}}^2 
+ | \boldsymbol{f}|_{L^2(0,T;\boldsymbol{H})}^2 \}$
and 
\begin{equation*}
	\nu \int_{0}^{t} 
	\bigl| \boldsymbol{v}(s) \bigr|_{\boldsymbol{V}}^2 ds
	\le M_8 \left( 1+ \frac{1}{\nu }
	\int_{0}^{T} \bigl| \tilde{\sigma }(s) \bigr|_{\mathbb{H}}^2 ds \right)
\end{equation*}
for all $t \in [0,T]$. Thus, we have \eqref{prov32} for some suitable constant $M_8>0$. 
\hfill $\Box$ \\

We can also define a solution operator $S_4:L^2(0,T;\mathbb{H})\to L^2(0,T;\boldsymbol{V})$ 
that assigns a unique solution $S_4 \tilde{\sigma}:=\boldsymbol{v}$ 
to the above problem \eqref{prov12}--\eqref{prov22} on $[0,T]$.

\subsection{Proof of Theorem 2.3.}

In this subsection, we prove Theorem~2.3 in a similar manner to the proof of Theorem~2.1.

\paragraph{Proof of Theorem~2.3.}

We define the operator $S:L^2(0,T;\boldsymbol{V}) \to L^2(0,T;\boldsymbol{V})$ 
by $S:=S_4\circ S_3$. 
Let $\nu \in (0,1]$. 
We show that $S$ is a contraction operator. 
For $\tilde{\boldsymbol{v}}^{(i)} \in L^2(0,T;\boldsymbol{V})$ for $i=1,2$, 
we set $\sigma^{(i)}:=S_3 \tilde{\boldsymbol{v}}^{(i)}$. 
Taking the difference between \eqref{prosig12} with $\sigma^{(1)}$ and 
{\risei \eqref{prosig12}} with $\sigma^{(2)}$ at {\fukao time $t=s$} and using 
the monotonicity of the subdifferential, we have 
\begin{align*}
	\frac{1}{2} \frac{d}{ds} \bigl| \sigma^{(1)}(s)-\sigma^{(2)}(s) \bigr|_{\mathbb{H}}^2 
	\le {} & \bigl( \varepsilon ( \tilde{\boldsymbol{v}}^{(1)}(s) - \tilde{\boldsymbol{v}}^{(2)}(s) \bigr),\sigma^{(1)}(s)-\sigma^{(2)}(s) \bigr) _{\mathbb{H}}
	\nonumber \\
	\le {} & \frac{1}{2} \bigl|  \tilde{\boldsymbol{v}}^{(1)}(s) - \tilde{\boldsymbol{v}}^{(2)}(s) \bigr|_{\boldsymbol{V}}^2
	+ \frac{1}{2} \bigl|  \sigma^{(1)}(s)-\sigma^{(2)}(s) \bigr|_{\mathbb{H}}^2
\end{align*} 
for a.a.\ $s \in (0,T)$. 
Using the {G}ronwall inequality, we obtain
\begin{equation}
	\bigl| \sigma^{(1)}(t)-\sigma^{(2)}(t) \bigr|_{\mathbb{H}}^2 
	\le  e^T
	\int_0^t
	\bigl| \tilde{\boldsymbol{v}}^{(1)}(s) 
	- \tilde{\boldsymbol{v}}^{(2)}(s) \bigr|_{\boldsymbol{V}}^2ds
	\label{banach5}
\end{equation}
for all $t \in [0,T]$. 

Next, we set $\boldsymbol{v}^{(i)}:=S_4 \sigma^{(i)}$ for 
$i=1,2$. 
In the same way, taking the difference between \eqref{prov12} with 
$\boldsymbol{v}^{(1)}$ and \eqref{prov12} with $\boldsymbol{v}^{(2)}$ at 
{\fukao time $t=s$}, we have
\begin{align*}
	& \frac{1}{2} \frac{d}{ds} \bigl| \boldsymbol{v}^{(1)}(s)-\boldsymbol{v}^{(2)}(s) \bigr|^2_{\boldsymbol{H}}
	+ \nu \bigl| \boldsymbol{v}^{(1)}(s) - \boldsymbol{v}^{(2)}(s) \bigr|_{\boldsymbol{V}}^2 \notag \\
	& \le - \bigl( \sigma ^{(1)}(s)-\sigma ^{(2)}(s), \varepsilon \bigl( \boldsymbol{v}^{(1)}(s)-\boldsymbol{v}^{(2)}(s) \bigr) \bigr)_{\mathbb{H}} \notag \\
	& \le \frac{1}{2\nu } \bigl| \sigma ^{(1)}(s)-\sigma ^{(2)}(s) \bigr|_{\mathbb{H}}^2 
	+ \frac{\nu }{2} \bigl| \boldsymbol{v}^{(1)}(s)-\boldsymbol{v}^{(2)}(s) \bigr|_{\boldsymbol{V}}^2,
\end{align*}
that is,
\begin{equation}
	\bigl| \boldsymbol{v}^{(1)}(t)-\boldsymbol{v}^{(2)}(t) \bigr|^2_{\boldsymbol{H}}
	+\nu \int_{0}^{t} \bigl| \boldsymbol{v}^{(1)}(s) - \boldsymbol{v}^{(2)}(s) \bigr|_{\boldsymbol{V}}^2 ds 
	\le \frac{1}{\nu } \int_{0}^{t} \bigl| \sigma ^{(1)}(s)-\sigma ^{(2)}(s) \bigr|_{\mathbb{H}}^2 ds
	\label{banach4}
\end{equation} 
for all $t \in [0,T]$. Combining \eqref{banach5} with {\risei \eqref{banach4}}, we obtain
the following inequality:
\begin{align*}
	\nu \bigl| S \tilde{\boldsymbol{v}}^{(1)}
	-S \tilde{\boldsymbol{v}}^{(2)}\bigr|^2_{L^2(0,T;\boldsymbol{V})} 
	\le {} & \frac{e^T}{\nu }T  \bigl| 
	\tilde{\boldsymbol{v}}^{(1)} - \tilde{\boldsymbol{v}}^{(2)} 
	\bigr|_{L^2(0,T;\boldsymbol{V})}^2.
\end{align*} 
Taking $T_0 \in (0,T]$ satisfying 
$(e^T /\nu^2)T_0 <1$, 
we see that $S$ is a contraction on $L^2(0,T_0;\boldsymbol{V})$.
Thus, we can apply the method of proof used in Theorem~2.1 to show Theorem~2.3. 
\hfill $\Box$

\paragraph{Remark.} 
The approach to the {\fukao viscous} perfect plasticity model is the same as in previous work
\cite[Chapter~V]{DL76}. 
There, the limiting problem $\nu \to 0$ was also considered. 
However, 
when the constraint set depends on time, 
our proof strategy does not work for 
the problem $\nu \to 0$. This is because 
we do not know how to treat the term 
\begin{equation*} 
	 \bigl({\risei \bigl( \partial I_{K(t)}\bigl( \sigma (t) \bigr)\bigr)'}, \sigma '(t) \bigr)_{\mathbb{H}}
\end{equation*} 
rigorously (cf.\ \cite[p.243]{DL76}). 
{\fukao Moreover, it is not possible to prove Theorem 2.3 under the assumption (A1)--(A5) by taking the 
limit $\kappa \to 0$ in Theorem 2.1. It is also not possible to prove 
Theorem 2.3 under the assumption (A1)--(A3), (A4$'$) and (A5) by taking the 
limit $\kappa \to 0$ in Theorem 2.2. }

\section*{Acknowledgments}

The present paper benefits from the support of the JSPS KAKENHI 
Grant-in-Aid for Scientific Research(C), Grant Number 17K05321 for TF.



\begin{thebibliography}{99}



\bibitem{Bir74}
	M.\ {B}iroli, 
	\newblock Sur les in\'equations paraboliques avec convexe d\'ependant du temps: 
	solution forte et solution faible, 
	\newblock Riv.\ Mat.\ Univ.\ Parma (3), \textbf{3} (1974), 33--72. 

\bibitem{Bir74b}
	M.\ {B}iroli, 
	\newblock Sur une inequation parabolique avec convexe dependant du temps, 
	\newblock Ricerche Mat., \textbf{23} (1974), 203--222.
	
\bibitem{BBW15} 
	S.\ {B}oettcher, M.\ {B}\"{o}hm and M.\ {W}olff, 
	\newblock Well-posedness of a thermo-elasto-plastic problem with phase transitions in TRIP steels under mixed boundary conditions, 
	\newblock ZAMM Z.\ Angew.\ Math.\ Mech., \textbf{95} (2015), 1461--1476. 

\bibitem{Bre72}
	H.\ {B}r\'ezis, 
	\newblock Un probl\'eme d'evolution avec contraintes unilat\'erales d\'ependant du temps, 
	\newblock C.\ R.\ Acad.\ Sci.\ Paris, S\'er A, \textbf{274} (1972), 310--312.

{\fukao 
\bibitem{Bre73}
	H.\ {B}r\'ezis, 
	\newblock {\it Op\'erateurs maximaux monotones et semi-groupes de contractions dans les especes de {H}ilbert}, 
	\newblock North-Holland, Amsterdam, 1973.
}

\bibitem{BKS04} 
	M.\ {B}rokate, P.\ {K}rej\v{c}\'{i} and H.\ {S}chnabel, 
	\newblock On uniqueness in evolution quasivariational inequalities, 
	\newblock J.\ Convex Anal., \textbf{11} (2004), 111--130. 

\bibitem{CR06}
	K.\ {C}he{\l}mi\'nski and R.\ {R}acke,
	\newblock Mathematical analysis of thermoplasticity with linear kinematic hardening, 
	\newblock J.\ Appl.\ Anal., \textbf{12} (2006), 37--57.

\bibitem{DL76} 
	G.\ {D}uvaut and J.\ L.\ {L}ions, 
	\newblock \textit{Inequalities in mechanics and physics}, 
	\newblock Springer-Verlag, 1976.

\bibitem{FK13}
	T.\ {F}ukao and N.\ {K}enmochi, 
	\newblock Abstract theory of variational inequalities and {L}agrange multipliers, 
	\newblock pp.\ 237--246 in \textit{Discrete and continuous dynamical systems, supplement 2013}, 2013.

\bibitem{FK13b}
	T.\ {F}ukao and N.\ {K}enmochi, 
	\newblock Parabolic variational inequalities with weakly time-dependent constraint, 
	\newblock Adv.\ Math.\ Sci.\ Appl., \textbf{23} (2013), 365--395.

\bibitem{KMK09} 
	R.\ {K}ano, Y.\ {M}urase and N.\ {K}enmochi, 
	\newblock Nonlinear evolution equations generated by subdifferentials with nonlocal constraints, 
	\newblock pp.175--194 in 
	{\it Nonlocal and abstract parabolic equations and their applications}, 
	{\bf Vol.86}, Banach Center Publications, 2009.

\bibitem{Ken75} 
	N.\ {K}enmochi, 
	\newblock Some nonlinear parabolic variational inequalities, 
	\newblock Israel J.\ Math., \textbf{22} (1975), 304--331. 

\bibitem{Ken81} 
	N.\ {K}enmochi,
	\newblock Solvability of nonlinear evolution equations with time-dependent constraints and applications, 
	\newblock Bull.\ Fac.\ Ed.\ Chiba Univ., \textbf{30} (1981), 1--87, 

{\fukao
\bibitem{Kis17}
	K.\ {K}isiel, 
	\newblock Dynamical poroplasticity model --- Existence theory for gradient type nonlinearities with {L}ipschitz perturbations, 
	\newblock J.\ Math.\ Anal.\ Appl., \textbf{450} (2017), 544--577. 
}

{\fukao
\bibitem{KK16}
	K.\ {K}isiel and K.\ {K}osiba, 
	\newblock Dynamical poroplasticity model with mixed boundary conditions --- Theory for ${\mathcal LM}$-type nonlinearity, 
	\newblock J.\ Math.\ Anal.\ Appl., \textbf{443} (2016), 187--229. 
}

\bibitem{Kre96} 
	P.\ {K}rej\v{c}\'{i}, 
	\newblock \textit{Hysteresis, convexity and dissipation in hyperbolic equations}, 
	\newblock Gakuto Intarnat.\ Ser.\ Math.\ Sci.\ Appl., \textbf{8} Gakkotosho, Tokyo, 1996.

\bibitem{KL02} 
	P.\ {K}rej\v{c}\'{i} and Ph.\ {L}auren\c{c}ot, 
	\newblock Generalized variational inequalities, 
	\newblock J.\ Convex Anal., \textbf{9} (2002), 159--183. 
	
\bibitem{KL09} 
	P.\ {K}rej\v{c}\'{i} and M.\ {L}iero, 
	\newblock Rate independent {K}urzweil processes, 
	\newblock Appl.\ Math., \textbf{54} (2009), 117--145. 

\bibitem{KR11} 
	P.\ {K}rej\v{c}\'{i} and T.\ {R}oche, 
	\newblock Lipschitz continuous data dependence of sweeping processes in BV spaces, 
	\newblock Discrete Contin.\ Dyn.\ Syst.\ Ser.\ B, \textbf{15} (2011), 637--650. 

\bibitem{KM98} 
	M.\ {K}unze and M.\ D.\ P.\ {M}onteiro {M}arques, 
	\newblock On parabolic quasi-variational inequalities and state-dependent sweeping processes, 
	\newblock Topol.\ Methods Nonlinear Anal., \textbf{12} (1998), 179--191.

\bibitem{Mon93} 
	M.\ D.\ P.\ {M}onteiro {M}arques, 
	\newblock \textit{Differential inclusions in nonsmooth mechanical problems, shocks and dry friction}, 
	\newblock Progr.\ Nonlinear Differential Equations Appl., Birkh{\"a}user, Boston, 1993.

\bibitem{Mor71} 
	J.\ J.\ {M}oreau, 
	\newblock Rafle par un convexe variable (Premi\'ere partie), 
	\newblock pp.\ 1--43 in \textit{Travaux du S\'eminaire d'Analyse Convexe}, Montpellier, 1971. 

\bibitem{Mor77} 
	J.\ J.\ {M}oreau, 
	\newblock Evolution problem associated with a moving convex set in a {H}ilbert spaces, 
	\newblock J.\ Differential Equations, \textbf{26} (1977), 347--374.

\bibitem{Rec11} 
	V.\ {R}ecupero, 
	\newblock A continuity method for sweeping processes, 
	\newblock J.\ Differential Equations, \textbf{251} (2011), 2125--2142.
	
\bibitem{Rec15} 
	V.\ {R}ecupero, 
	\newblock BV continuous sweeping processes, 
	\newblock J.\ Differential Equations, \textbf{259} (2015), 4253--4272.

\bibitem{Sim87}
	J.\ {S}imon, 
	\newblock Compact sets in the spaces $L^p(0,T;B)$, 
	\newblock Ann.\ Mat.\ Pura.\ Appl.~(4), {\bf 146} (1987), 65--96.

\bibitem{Sho97}
	R.\ E.\ {S}howalter, 
	\newblock \textit{Monotone operators in {B}anach spaces and nonlinear partial differential equations},
	\newblock Mathematical Surveys and Monographs, Vol.49, American Mathematical Society, 1997.

\bibitem{Ste04} 
	U.\ {S}tefanelli, 
	\newblock Some quasivariational problems with memory, 
	\newblock Boll.\ Unione Mat.\ Ital.\ Sez.\ B Artic.\ Ric.\ Mat.\ (8), \textbf{7} (2004), 319--333.

\bibitem{Vla13}
	A.\ {V}ladimirov,
	\newblock Equicontinuous sweeping processes, 
	\newblock Discrete Contin.\ Dyn.\ Syst.\ Ser.\ B, \textbf{18} (2013), 565--573. 

\bibitem{Yam76} 
	Y.\ {Y}amada, 
	\newblock On evolution equations generated by subdifferential operators, 
	\newblock J.\ Fac.\ Sci.\ Univ.\ Tokyo Sect.\ IA Math., {\bf 23} (1976), 491--515.

\end{thebibliography}
\end{document}